\documentclass[review,onefignum,onetabnum]{siamart171218}

\usepackage{lipsum}
\usepackage{mathrsfs}
\usepackage{amsfonts}
\usepackage{graphicx}
\usepackage{epstopdf}
\usepackage{algorithmic}
\usepackage{amssymb}
\usepackage{amsmath}
\usepackage{adjustbox}

\DeclareMathAlphabet{\mathpzc}{OT1}{pzc}{m}{it}

\ifpdf
  \DeclareGraphicsExtensions{.eps,.pdf,.png,.jpg}
\else
  \DeclareGraphicsExtensions{.eps}
\fi

\newsiamremark{hypothesis}{Hypothesis}
\crefname{hypothesis}{Hypothesis}{Hypotheses}
\newsiamthm{claim}{Claim}

\nolinenumbers

\headers{Zeros of parabolic cylinder functions}{T. M. Dunster, A. Gil, D. Ruiz-Antolin and J. Segura}

\title{Uniform asymptotic expansions for the zeros of parabolic cylinder functions}

\author{T. M. Dunster\thanks{Department of Mathematics and Statistics, San Diego State University, 5500 Campanile Drive, San Diego, CA 92182, USA. 
  (\email{mdunster@sdsu.edu}, \url{https://tmdunster.sdsu.edu}).}
\and A. Gil\thanks{Departamento de Matem\'atica Aplicada y CC. de la Computaci\'on, ETSI Caminos, Universidad de Cantabria, 39005-Santander, Spain. 
  (\email{amparo.gil@unican.es}).}
  \and D. Ruiz-Antolin\thanks{Departamento de Matem\'atica Aplicada y CC. de la Computaci\'on, ETSI Caminos, Universidad de Cantabria, 39005-Santander, Spain. 
  (\email{diego.ruizantolin@unican.es}).}
\and J. Segura\thanks{Departamento de Matem\'aticas, Estad\'{\i}stica y Computaci\'on, Facultad de Ciencias, Universidad de Cantabria, 39005-Santander, Spain. 
  (\email{javier.segura@unican.es}).}}
\usepackage{amsopn}

\makeatletter
\newcommand*{\addFileDependency}[1]{% argument=file name and extension
  \typeout{(#1)}% latexmk will find this if $recorder=0 (however, in that case, it will ignore #1 if it is a .aux or .pdf file etc and it exists! if it doesn't exist, it will appear in the list of dependents regardless)
  \@addtofilelist{#1}% if you want it to appear in \listfiles, not really necessary and latexmk doesn't use this
  \IfFileExists{#1}{}{\typeout{No file #1.}}% latexmk will find this message if #1 doesn't exist (yet)
}
\makeatother

\nolinenumbers

\ifpdf
\hypersetup{
  pdftitle={Uniform asymptotic expansions for the zeros of parabolic cylinder functions},
  pdfauthor={T. M. Dunster, A. Gil, D. Ruiz-Antolin and J. Segura}
}

\begin{document}

\maketitle

\begin{abstract}
The real and complex zeros of the parabolic cylinder function $U(a,z)$ are studied. Asymptotic expansions for the zeros are derived, involving the zeros of Airy functions, and these are valid for $a$ positive or negative and large in absolute value, uniformly for unbounded $z$ (real or complex). The accuracy of the approximations of the complex zeros is then demonstrated with some comparative tests using a highly precise numerical algorithm for finding the complex zeros of the function.
\end{abstract}



\section{Introduction} 
\label{eqsec1}

We consider the zeros of the parabolic cylinder function \allowbreak 
$U(a,z)$, where $a \in \mathbb{R}$ and $z \in \mathbb{C}$. It is a solution of differential equation
%%%%%%%%%%%%%%%%%%%
\begin{equation}  
\label{eq01}
\frac{d^{2}y}{dz^{2}}-\left(\frac{1}{4}z^{2}+a \right)y=0,
\end{equation}
%%%%%%%%%%%%%%%%%%%
and has the integral representations \cite[Eqs. 12.5.1 and 12.5.4]{NIST:DLMF}
%%%%%%%%%%%%%%%%%%%
\begin{equation}
\label{eq02}
U(a,z)=\frac{e^{-\frac{1}{4}z^{2}}}{\Gamma
\left(a+\frac{1}{2} \right)}
\int_0^\infty {t^{a-\frac{1}{2}}
e^{-\frac{1}{2}t^{2}-zt}dt}
\quad \left(a>-\tfrac{1}{2} \right),
\end{equation}
%%%%%%%%%%%%%%%%%%%
and
%%%%%%%%%%%%%%%%%%%
\begin{equation}
\label{eq03}
U(-a,z)=\sqrt {\frac{2}{\pi }} e^{
\frac{1}{4}z^{2}}\int_0^\infty {t^{a-\frac{1}{2}}e^{-\frac{1}{2}t^{2}}\cos \left( {zt-\tfrac{1}{2}\pi a+\tfrac{1}{4}\pi } \right)dt} \quad \left(a>-\tfrac{1}{2} \right).
\end{equation}
%%%%%%%%%%%%%%%%%%%
It can also be expressed in terms of the confluent hypergeometric function \cite[Sect. 13.2]{NIST:DLMF} by
%%%%%%%%%%%%%%%%%%%
\begin{equation}
\label{eq04}
U(a,z)=2^{-\frac{1}{4}(2a+1)}
e^{-\frac{1}{4}z^{2}} U\left( \tfrac{1}{2}a+\tfrac{1}{4},
\tfrac{1}{2},\tfrac{1}{2}z^{2} \right),
\end{equation}
%%%%%%%%%%%%%%%%%%%
and for the special case $a=-n-\frac12$ ($n=0,1,2,\ldots$) the function is a multiple of the Hermite polynomial of order $n$ (see \cite[Eq. 12.7.2 and Sec. 18.3.1]{NIST:DLMF})
%%%%%%%%%%%%%%%%%%%
\begin{equation}
\label{eq04a}
U(-n-\tfrac12,z)=e^{-\frac14 z^2}He_{n}(z)
=2^{-\frac12 n} e^{-\frac14 z^2} H_{n}(2^{-\frac12}z).
\end{equation}
%%%%%%%%%%%%%%%%%%%

The important property of $U(a,z)$ is that it is the unique solution of \eqref{eq01} which is recessive in the sector $|\arg(z)|\leq \pi /4$, with all other independent solutions being dominant in this sector. This follows from
%%%%%%%%%%%%%
\begin{equation}
\label{eq05}
U(a,z)\sim z^{-a-\frac{1}{2}}
e^{-\frac{1}{4}z^{2}}
\quad \left(z\to \infty, \, \vert \arg (z)\vert \le
\tfrac{3}{4}\pi -\delta 
\right).
\end{equation}
%%%%%%%%%%%%%%%%%%%

The zeros of parabolic cylinder functions find applications in different areas of science and engineering. For example, the complex zeros 
of parabolic cylinder functions are used in the analysis of fluid flow in cylindrical channels, such as in the study of laminar flow or boundary layer theory
or in the scattering of electromagnetic waves at parabolic boundaries \cite{Hillion:1997:EWP}.

In this paper we give asymptotic expansions for the zeros of $U(a,z)$ for $a$ positive (\cref{sec:apos}) and $a$ negative (\cref{sec:aneg}). These are valid for $|a|$ large, uniformly for unbounded $z$ (real or complex). Previous uniform approximations for the real zeros of $U(a,x)$ with $|a|$ large were given by Olver \cite{Olver:1959:USE}, but his expansions are difficult to compute beyond two terms. Here we extend Olver's results by obtaining full expansions with coefficients that are explicitly given and readily computable, and we also consider the complex zeros. We use the method that was recently developed in \cite{Dunster:2024:AZB} for the uniform asymptotic expansion of the zeros of Bessel functions.  

In \cref{num:tests} the accuracy of the approximations to the complex zeros is tested using a numerical implementation of a method for finding the complex zeros of solutions of second order ODEs described in \cite{Segura:2013:CCZ}.

\section{Zeros of \texorpdfstring{$U(a,z)$}{} for \texorpdfstring{$a>0$}{}}
\label{sec:apos}

Following \cite[Eq. (3.3)]{Dunster:2021:UAP} let
%%%%%%%%%%
\begin{equation}
\label{eqZeta}
\frac{2}{3}\zeta^{3/2}=\xi
=\int_{1}^{\hat{z}} \left(t^2-1\right)^{1/2} dt
=\frac{1}{2}\hat{z}\left(\hat{z}^{2}-1\right)^{1/2}
-\frac{1}{2}\ln\left\{\hat{z}+\left(\hat{z}^{2}-1\right)^{1/2} \right\},
\end{equation}
%%%%%%%%%%
where principal branches apply. Here $\xi$ and $\zeta$ are the variables associated with Liouville-Green and turning point problems, respectively (see \cite[Chaps. 10 and 11]{Olver:1997:ASF}). $\zeta$ is an analytic function of $\hat{z}$ in the plane having a branch cut along the interval $(-\infty,-1]$, and $\xi$ is an analytic function of $\hat{z}$ in the plane having a branch cut along the interval $(-\infty,1]$. Moreover $\xi,\, \zeta \in [0,\infty)$ for $1\leq \hat{z}<\infty$. Also $\zeta \leq 0$ for $-1\leq \hat{z}\leq 1$ with
%%%%%%%%%%
\begin{equation} 
\label{eqZeta1}
\tfrac{2}{3}(-\zeta)^{3/2}
=\tfrac{1}{2}\arccos(\hat{z})
-\tfrac{1}{2}\hat{z}\left(1-\hat{z}^{2}\right)^{1/2}.
\end{equation}
%%%%%%%%%%
Note $\hat{z} = 0$ corresponds to $\zeta = -\frac14(3\pi)^{2/3}$.

Now from \cite[Eqs. (3.15) - (3.17), (3.35) and (3.36)]{Dunster:2021:UAP}
%%%%%%%%%%
\begin{multline} 
\label{eq3.29}
U\left(-\tfrac{1}{2}u,\sqrt {2u}\,\hat{z}\right)
= \sqrt {\pi }\,2^{\frac{1}{4}(3-u)}
e^{-\frac{1}{4}u} u^{\frac{1}{12}(3u-1)}
\left(\frac{\zeta}{\hat{z}^2-1}\right)^{1/4}
\\
\times 
\left[\left\{1+\mathcal{O}\left( \frac{1}{u^2} \right)\right\}
\mathrm{Ai}\left( u^{2/3}\zeta \right)
+\mathcal{O}\left( \frac{1}{u^{4/3} (|\zeta|+1)^{1/2}} \right)
\mathrm{Ai}'\left( u^{2/3}\zeta \right)
\right],
\end{multline}
%%%%%%%%%%
and
%%%%%%%%%%
\begin{multline} 
\label{eq3.29a}
U\left(\tfrac{1}{2}u,\pm i\sqrt {2u}\,\hat{z}\right)
= \sqrt {\pi } \, 2^{\frac{1}{4}(3+u)}
u^{-\frac{1}{12}(3u+1)} 
\exp\left\{\frac{1}{4}(1 \mp \pi i)u 
\mp\frac{1}{12}\pi i\right\} 
\left(\frac{\zeta}{\hat{z}^2-1}\right)^{1/4}
\\
\times 
\left[\left\{1+\mathcal{O}\left( \frac{1}{u^2} \right)\right\}
\mathrm{Ai}_{\mp 1}\left( u^{2/3}\zeta \right)
+\mathcal{O}\left( \frac{1}{u^{4/3} (|\zeta|+1)^{1/2}} \right)
\mathrm{Ai}_{\mp 1}'\left( u^{2/3}\zeta \right)
\right],
\end{multline}
%%%%%%%%%%
uniformly in a certain complex domain which contains the half plane $\Re(\hat{z}) \geq 0$. Here we use the standard notation $\mathrm{Ai}_{l}(z):=\mathrm{Ai}(ze^{-2\pi il/3})$ ($l=\pm 1$).

Based on (\ref{eq3.29a}), and following \cite{Dunster:2024:AZB}, we define $\mathcal{Y}(u,\hat{z})$ and $\mathcal{Z}(u,\hat{z})$ by the pair of equations
%%%%%%%%%%
\begin{equation}
\label{eq05a}
U\left(\tfrac{1}{2}u, \pm i\sqrt {2u} \, \hat{z}\right)
=e^{\mp \frac{1}{12}(3u+1)\pi i} \mathcal{Y}(u,\hat{z})
\mathrm{Ai}_{\mp 1}\left( u^{2/3}\mathcal{Z}(u,\hat{z}) \right).
\end{equation}
%%%%%%%%%%

To see why the factors $e^{\mp (3u+1)\pi i/12}$ are required, we use the Airy function connection formula \cite[Eq. 9.2.12]{NIST:DLMF}
%%%%%%%%%%%%%%%%%%%
\begin{equation}
\label{eqAiConnect}
\mathrm{Ai}(z)
=e^{\pi i/3}\mathrm{Ai}_{1}(z)
+e^{-\pi i/3}\mathrm{Ai}_{-1}(z),
\end{equation}
%%%%%%%%%%%%%%%%%%%
and the parabolic cylinder function connection formula \cite[Eq. 12.2.18]{NIST:DLMF}
%%%%%%%%%%
\begin{multline}
\label{eq1.8}
U\left(-\tfrac{1}{2}u,\pm \sqrt {2u} \, \hat{z}\right)
=\frac{\Gamma\left(\tfrac{1}{2}u+\tfrac{1}{2}\right)}{\sqrt{2\pi}}
\\
\times 
\left\{e^{ \frac{1}{4}(u-1)\pi i}
U\left(\tfrac{1}{2}u,\pm i \sqrt {2u} \, \hat{z}\right)
+ e^{- \frac{1}{4}(u-1)\pi i}
U\left(\tfrac{1}{2}u,\mp i \sqrt {2u} \, \hat{z}\right)\right\}.
\end{multline}
%%%%%%%%%%
Then from (\ref{eq05a}) - (\ref{eq1.8}) 
%%%%%%%%%%%%%%%%%%%
\begin{equation}
\label{eq1.8a}
U\left(-\tfrac{1}{2}u,\sqrt {2u} \, \hat{z}\right)
=\frac{\Gamma\left(\tfrac{1}{2}u+\tfrac{1}{2}\right)}{\sqrt{2\pi}}\mathcal{Y}(u,\hat{z})
\mathrm{Ai}\left( u^{2/3}
\mathcal{Z}(u,\hat{z}) \right),
\end{equation}
%%%%%%%%%%%%%%%%%%%
which matches the form of the Airy approximation (\ref{eq3.29}). From \cite{Dunster:2024:AZB} it is shown that $\mathcal{Y}(u,\hat{z})$ and $\mathcal{Z}(u,\hat{z})$ are bounded analytic functions in the intersection of the domains for which (\ref{eq3.29}) and (\ref{eq3.29a}) are valid, which includes the half-plane $\Re(\hat{z}) \geq 0$. 

To find $\mathcal{Y}(u,\hat{z})$ use (\ref{eq05a}) to get
%%%%%%%%%%
\begin{multline}
\label{eq1.9}
i\sqrt {2u}\left\{
U\left(\tfrac{1}{2}u,  i\sqrt {2u} \, \hat{z}\right)
U'\left(\tfrac{1}{2}u, -i\sqrt {2u} \, \hat{z}\right)
+U'\left(\tfrac{1}{2}u, i\sqrt {2u} \, \hat{z}\right)
U\left(\tfrac{1}{2}u,- i\sqrt {2u} \, \hat{z}\right)\right\}
\\
=  u^{2/3} \mathcal{Y}^2(u,\hat{z})\mathcal{Z}'(u,\hat{z})
\left\{
\mathrm{Ai}_{1}\left( u^{2/3}\mathcal{Z}(u,\hat{z}) \right)
\mathrm{Ai}'_{-1}\left( u^{2/3}\mathcal{Z}(u,\hat{z}) \right)
\right.
\\  \left.
-\mathrm{Ai}'_{1}\left( u^{2/3}\mathcal{Z}(u,\hat{z}) \right)
\mathrm{Ai}_{-1}\left( u^{2/3}\mathcal{Z}(u,\hat{z}) \right)
\right\},
\end{multline}
%%%%%%%%%%
where $\mathcal{Z}'(u,\hat{z}) =\partial  \mathcal{Z}(u,\hat{z})/ \partial \hat{z}$. Then use the Wronskian relations \cite[Eqs. 9.2.9 and 12.2.11]{NIST:DLMF} and we arrive at
%%%%%%%%%%%%%%%%%%%
\begin{equation}
\label{eq1.10}
\mathcal{Y}(u,\hat{z})
=2 \pi^{3/4} u^{-1/12}
\left\{\Gamma\left(\tfrac{1}{2}u+\tfrac{1}{2}\right)
\mathcal{Z}'(u,\hat{z})\right \}^{-1/2}.
\end{equation}
%%%%%%%%%%%%%%%%%%%

The function $\mathcal{Z}(u,\hat{z})$ is more important in the study of the zeros, and from \cite[Eq. (2.14)]{Dunster:2024:AZB}
%%%%%%%%%%%%%%%%%%%
\begin{equation}
\label{eq1.11}
\mathcal{Z}(u,\hat{z}) \sim \zeta
+\sum_{s=1}^{\infty}
\frac{\Upsilon_{s}(\hat{z})}{u^{2s}}
\quad  (u \to \infty),
\end{equation}
%%%%%%%%%%%%%%%%%%%
uniformly (in our case) for an unbounded domain containing $\Re(\hat{z}) \geq 0$, with coefficients $\Upsilon_{s}(\hat{z})$ described next.

In order to do so, we first must give formulas for certain Liouville-Green coefficients which appear in \cite[Eqs. (2.4) and (2.5)]{Dunster:2024:AZB}. To this end define
%%%%%%%%%%
\begin{equation}
\label{eq06}
\beta=\hat{z}\left(\hat{z}^{2}-1\right)^{-1/2},
\end{equation}
%%%%%%%%%%
where the principal branch of the square root is taken, so that $\beta$ is positive for $\hat{z}>1$ and is continuous in the plane having a cut along $[-1,1]$. Thus $\beta \rightarrow 1$ as $\hat{z} \rightarrow \infty$ in any direction. Next from \cite[Eqs. (3.7) - (3.9)]{Dunster:2021:UAP} define
%%%%%%%%%%%%%
\begin{equation}
\label{eq07}
\mathrm{E}_{1}(\beta)=\tfrac{1}{24}\beta
\left(5\beta^{2}-6\right),
\end{equation}
%%%%%%%%%%%%%%%%
\begin{equation}
\label{eq08}
\mathrm{E}_{2}(\beta)=
\tfrac{1}{16}\left(\beta^{2}-1\right)^{2} 
\left(5\beta^{2}-2\right),
\end{equation}
and for $s=2,3,4\cdots$
%%%%%%%%%%%%%%%%
\begin{equation}
\label{eq09}
\mathrm{E}_{s+1}(\beta) =
\frac{1}{2} \left(\beta^{2}-1 \right)^{2}\mathrm{E}_{s}^{\prime}(\beta)
+\frac{1}{2}\int_{\omega(s)}^{\beta}
\left(p^{2}-1 \right)^{2}
\sum\limits_{j=1}^{s-1}
\mathrm{E}_{j}^{\prime}(p)
\mathrm{E}_{s-j}^{\prime}(p) dp,
\end{equation}
%%%%%%%%%%%%%%%%
where $\omega(s)=1$ for $s$ odd and $\omega(s)=0$ for $s$ even. We remark that $\mathrm{E}_{2s}(\beta)$ is even, $\mathrm{E}_{2s+1}(\beta)$ is odd, and $\mathrm{E}_{2s}(\pm 1)=0$. Then the coefficients used in \cite[Eqs. (2.4) and (2.5)]{Dunster:2024:AZB} are $\hat{E}_{s}(\hat{z})=\mathrm{E}_{s}(\beta)$ (here we use $\hat{z}$ in place of $z$).

We only need the odd coefficients. Let 
%%%%%%%%%
\begin{equation}
\label{eqG}
G_{s}(\hat{z})=\beta^{-1}E_{2s-1}(\hat{z})
\quad  (s=1,2,3,\ldots),
\end{equation}
%%%%%%%%%
which are even polynomials in terms of $\beta$, and hence from (\ref{eq06}) rational functions of $\hat{z}$. The first four are found to be
%%%%%%%%%
\begin{equation}
\label{eqG1}
G_{1}(\hat{z})
=-\frac{\hat{z}^2-6}{24\left(\hat{z}^2-1\right)},
\end{equation}
%%%%%%%%%
%%%%%%%%%
\begin{equation}
\label{eqG2}
G_{2}(\hat{z})
=\frac{56 \hat{z}^8 - 252 \hat{z}^6 
+ 441 \hat{z}^4 + 1860 \hat{z}^2 
+ 3420}{5760\left(\hat{z}^2-1\right)^4},
\end{equation}
%%%%%%%%%
%%%%%%%%%
\begin{multline}
\label{eqG3}
G_{3}(\hat{z})
=-\frac{1}{322560\left(\hat{z}^2-1\right)^7}
\left(3968 \hat{z}^{14} - 29760 \hat{z}^{12} + 96720 \hat{z}^{10}
- 177320 \hat{z}^8  \right.
\\  \left.
+ 199485 \hat{z}^6 
- 1719018 \hat{z}^4 - 5480580 \hat{z}^2 - 1590120\right),
\end{multline}
%%%%%%%%%
and
%%%%%%%%%
\begin{multline}
\label{eqG4}
G_{4}(\hat{z})
=\frac{1}{3440640\left(\hat{z}^2-1\right)^{10}}
\left(130048 \hat{z}^{20} - 1365504\hat{z}^{18} 
+ 6486144 \hat{z}^{16} \right.
\\
- 18377408 \hat{z}^{14} 
+ 34457640 \hat{z}^{12} - 44794932 \hat{z}^{10} 
+ 41062021 \hat{z}^8 + 495103464 \hat{z}^6  
\\
\left.
+ 3107060712 \hat{z}^4 + 2497542880 \hat{z}^2 + 292852560 
\right).
\end{multline}
%%%%%%%%%

Now we can apply \cite[Thm. 2.2]{Dunster:2024:AZB} to construct the aforementioned coefficients $\Upsilon_{s}(\hat{z})$ of (\ref{eq1.11}). The first four are given explicitly by \cite[Eqs. (2.28) and (2.40) - (2.42)]{Dunster:2024:AZB}. In these we can use (\ref{eq06}) and (\ref{eqG}), and recalling $\xi =\frac{2}{3}\zeta^{3/2}$, to replace $\frac {3}{2} \xi \zeta^{-2} E_{2s+1}$ by $\hat{z} \sigma \zeta^{-1} G_{s+1}$ where
%%%%%%%%%%%
\begin{equation}
\label{eq20}
\sigma(\hat{z})
=\left( \frac{\zeta }{\hat{z}^2-1 }\right)^{1/2}.
\end{equation}
%%%%%%%%%%%
Thus for example
%%%%%%%%%%%
\begin{equation}
\label{eq20a}
\Upsilon_{1}(\hat{z})=
\frac{\hat{z} \sigma(\hat{z})  G_{1}(\hat{z})}{\zeta}
-\frac{5}{48\zeta^{2}}.
\end{equation}
%%%%%%%%%%%
The coefficients $\Upsilon_{s}(\hat{z})$ are analytic for $\Re(\hat{z}) \geq 0$, with $\hat{z}=1$ ($\zeta=0$) being removable singularities (as is true for $\sigma(\hat{z}))$.

Consider now the parabolic cylinder function with the upper signs in (\ref{eq05a}) (we consider the negative parameter case (\ref{eq3.29}) in the next section). We solve for the zeros $\hat{z} = \hat{z}_{m} \in \mathbb{C}$ ($m=1,2,3,\ldots$) asymptotically so that
%%%%%%%%%%
\begin{equation}
\label{eq21}
U\left(\tfrac{1}{2}u, i\sqrt {2u} \, \hat{z}_{m}\right)=0.
\end{equation}
%%%%%%%%%%

It is well known that the zeros of $U(a,z)$ lie in the second and third quadrants when $a > -\frac12$, and are conjugates of one another. Thus from (\ref{eq21}) it suffices for us to consider $\hat{z}_{m}$ lying in the first quadrant.

Now with $\mathrm{a}_{m}$ denoting the $m$th (negative) zero of $\mathrm{Ai}(x)$, ordered in increasing absolute values, we have from (\ref{eq05a}) and (\ref{eq21})
%%%%%%%%%%
\begin{equation}
\label{eq22}
e^{2\pi i/3} u^{2/3}\mathcal{Z}(u,\hat{z}_{m})
=\mathrm{a}_{m} 
\quad (m=1,2,3,\ldots),
\end{equation}
%%%%%%%%%%
recalling  $\mathrm{Ai}_{-1}(z)=\mathrm{Ai}(ze^{2\pi i/3})$.

From (\ref{eq1.11}) and (\ref{eq22}) we observe first order approximation is $\hat{z}_{m} \approx  \hat{z}_{m,0}$, where $\hat{z}_{m,0}$ satisfies
%%%%%%%%%%%%%%%%%%%
\begin{equation}
\label{eq33}
e^{2\pi i/3} u^{2/3} \zeta\left(\hat{z}_{m,0}\right)
=\mathrm{a}_{m} < 0,
\end{equation}
%%%%%%%%%%%%%%%%%%%
so that $\hat{z}_{m,0}$ is the solution lying in the first quadrant of the implicit equation
%%%%%%%%%%
\begin{multline}
\label{eq34}
2\int_{1}^{\hat{z}_{m,0}} \left(t^2-1\right)^{1/2} dt
=\hat{z}_{m,0}\left(\hat{z}_{m,0}^{2}-1\right)^{1/2}
-\ln\left\{\hat{z}_{m,0}
+\left(\hat{z}_{m,0}^{2}-1\right)^{1/2} \right\}
\\
=\frac{4}{3}\left\{
\zeta\left(\hat{z}_{m,0}\right)\right\}^{3/2}
=\frac{4i}{3u}\left|\mathrm{a}_{m}\right|^{3/2}.
\end{multline}
%%%%%%%%%%

For fixed $m$ and $u \to \infty$ we have $\hat{z}_{m} \to 1$, and more precisely
%%%%%%%%%%
\begin{equation}
%\label{eq}
\hat{z}_{m,0}=1+\frac{e^{\pi i/3}
\left|\mathrm{a}_{m}\right|}
{2^{1/3}u^{2/3}}
+\mathcal{O}\left(\frac{1}{u^{4/3}}\right).
\end{equation}
%%%%%%%%%%
As $U_{m}:=|\mathrm{a}_{m}|^{3/2}/u \to \infty$ we have $\hat{z}_{m,0} \to \infty$ such that
%%%%%%%%%%
\begin{equation}
%\label{eq}
\hat{z}_{m,0}=2 e^{\pi i/4}
\sqrt{\frac{U_{m}}{3}}
+\frac{e^{-\pi i/4} \sqrt{3}
\ln\left(U_{m}\right)}
{8 \sqrt{U_{m}}}
+\mathcal{O}\left(
\frac{1}{\sqrt{U_{m}}}\right).
\end{equation}
%%%%%%%%%%
Our results are uniformly valid for all $m$, and so are valid for $\hat{z}_{m}$ close to the turning point $\hat{z}=1$ as well as for unbounded values.

We denote the corresponding leading approximation to $\zeta$ by $\zeta_{m,0}$ (which lies in the first quadrant), and so from (\ref{eq33})
%%%%%%%%%%%%%%%%%%%
\begin{equation}
\label{eq35}
\zeta_{m,0}
=\zeta\left(\hat{z}_{m,0}\right)
=e^{\pi i/3} u^{-2/3}\left|\mathrm{a}_{m}\right|,
\end{equation}
%%%%%%%%%%%%%%%%%%%, 
and similarly define
%%%%%%%%%%%%%%%%%%%
\begin{equation}
\label{eq36a}
\zeta_{m,0}^{\prime}
=\zeta'\left(\hat{z}_{m,0}\right),  \;
\zeta_{m,0}^{\prime \prime}
=\zeta''\left(\hat{z}_{m,0}\right),  \;
\ldots.
\end{equation}
%%%%%%%%%%%%%%%%%%%
Primes here and elsewhere are derivatives with respect to $\hat{z}$.

Now the exact solution of (\ref{eq21}) has the asymptotic expansion
%%%%%%%%%%%%%%%%%%%
\begin{equation}
\label{eq50}
\hat{z}_{m} \sim \hat{z}_{m,0}+\sum_{s=1}^{\infty}
\frac{\hat{z}_{m,s}}{u^{2s}}
\quad (u \to \infty).
\end{equation}
%%%%%%%%%%%%%%%%%%%
The steps in determining the coefficients $z_{m,s}$ ($s=1,2,3,\ldots$) follow exactly as in those for Bessel functions \cite[Sect. 3]{Dunster:2024:AZB}, and the first four are given by the same general equations (3.47) - (3.50) of that paper, with $z$ replaced by $\hat{z}$. In these, for $s=1,2,3,\ldots$
%%%%%%%%%%%%%%%%%%%
\begin{equation}
\label{eq38a}
\Upsilon_{m,s}
=\Upsilon_{s}\left(\hat{z}_{m,0}\right),  \;
\Upsilon_{m,s}^{\prime}
=\Upsilon_{s}'\left(\hat{z}_{m,0}\right),  \;
\Upsilon_{m,s}^{\prime \prime}
=\Upsilon_{s}''\left(\hat{z}_{m,0}\right),  \;
\ldots.
\end{equation}
%%%%%%%%%%%%%%%%%%%
These of course differ than those for corresponding ones in the Bessel function case, since here $\hat{z}_{m,0}$ and $\Upsilon_{s}(\hat{z})$ are different.

To aid in the evaluation the derivatives, we have from (\ref{eqZeta})
%%%%%%%%%%%%%%%%%%%
\begin{equation}
\label{eq28}
\zeta'=
\left(\frac{\hat{z}^2-1}{\zeta}\right)^{1/2}
=\frac{1}{\sigma},
\end{equation}
%%%%%%%%%%%%%%%%%%%
and from (\ref{eq20})
%%%%%%%%%%%%%%%%%%%
\begin{equation}
\label{eqSigmap}
\sigma'=
\frac{1}{2 \zeta}
-\frac{\hat{z} \sigma}
{\hat{z}^{2}-1}
=\frac{1}{2 \zeta}
\left(1-2\hat{z} \sigma^{3}\right).
\end{equation}
%%%%%%%%%%%%%%%%%%%
Thus by the chain rule and repeated use of (\ref{eq28}) and (\ref{eqSigmap}), the derivatives of $\zeta(\hat{z})$ and $\sigma(\hat{z})$ at $\hat{z}=\hat{z}_{m,0}$ can all be expressed as rational functions of the computed values $\hat{z}_{m,0}$, $\zeta_{m,0}$ and $\sigma_{m,0}$, where
%%%%%%%%%%%%%%%%%%%
\begin{equation}
\label{eqSigma0}
\sigma_{m,0}=\sigma\left(z_{m,0}\right)
=\left( \frac{\zeta_{m,0} }
{\hat{z}_{m,0}^2-1 }\right)^{1/2}
=
\left(\frac{e^{\pi i/3}|\mathrm{a}_{m}|}
{ u^{2/3}\left(\hat{z}_{m,0}^{2}-1\right)}\right)^{1/2}.
\end{equation}
%%%%%%%%%%%%%%%%%%%

For example
%%%%%%%%%%%%%%%%%%%
\begin{equation}
%\label{eq}
\zeta_{m,0}''=
\frac{2\hat{z}_{m,0}\sigma_{m,0}^{3}-1}
{2\sigma_{m,0}^{2} \zeta_{m,0}},
\end{equation}
%%%%%%%%%%%%%%%%%%%
and
%%%%%%%%%%%%%%%%%%%
\begin{equation}
%\label{eq}
\sigma_{m,0}''=
\frac{6\sigma_{m,0}^{6}+4
\sigma_{m,0}^{4} \zeta_{m,0}
-\hat{z}_{m,0}\sigma_{m,0}^{3}-1}{2
\sigma_{m,0} \zeta_{m,0}^{2}}.
\end{equation}
%%%%%%%%%%%%%%%%%%
Using this method the first and second terms of in the series (\ref{eq50}) are found to be
%%%%%%%%%
\begin{equation}
\label{eqhatz1}
\hat{z}_{m,1}  =
\frac{\sigma_{m,0}}
{48 \zeta_{m,0}^{2}}
\left\{
12 \hat{z}_{m,0} \sigma_{m,0} \zeta_{m,0}
-10 \hat{z}_{m,0}^{3} \sigma_{m,0}^{3}
+5\right\},
\end{equation}
%%%%%%%%%%%
and
%%%%%%%%%%%%%%%%%%%
\begin{multline}
\label{eqhatz2}
\hat{z}_{m,2}=
-\frac{\sigma_{m,0}}{46080
\zeta_{m,0}^{5}}
\left\{
200 z_{m,0}^{7} \sigma_{m,0}^{9}
\left(221 z_{m,0}^{2}+35\right) 
\right.
\\
-720 z_{m,0}^{5}\sigma_{m,0}^{7}  \zeta_{m,0}
\left(221z_{m,0}^{2}+25\right)
-4000z_{m,0}^{4}\sigma_{m,0}^{6}
+ 24 z_{m,0}^{3}\sigma_{m,0}^{5}\zeta_{m,0}^{2}
\left(8847 z_{m,0}^{2} + 580\right)
\\
+5400 z_{m,0}^{2} \sigma_{m,0}^{4}\zeta_{m,0}
-10z_{m,0} \sigma_{m,0}^{3} 
\left(12432 z_{m,0}^{2} \zeta_{m,0}^{3}
+ 288 \zeta_{m,0}^{3} - 25\right)
 \\
 \left. 
 - 1200 \sigma_{m,0}^{2} \zeta_{m,0}^{2}
+27360  z_{m,0} \sigma_{m,0} \zeta_{m,0}^{4} 
-5525 \right\},
\end{multline}
%%%%%%%%%%%%%%%%%%%
and so on. Compare these with the corresponding coefficients for Bessel functions \cite[Eqs. (3.56) and (3.57)]{Dunster:2024:AZB}.

\section{Zeros of \texorpdfstring{$U(a,z)$}{} for \texorpdfstring{$a<0$}{}}
\label{sec:aneg}

Consider next
%%%%%%%%%%
\begin{equation}
\label{eq50a}
U\left(-\tfrac{1}{2}u,  \sqrt {2u} \, \hat{z}\right)=0,
\end{equation}
%%%%%%%%%%
where $u>0$. From \cite[Sect. 12.11(i)]{NIST:DLMF} there are no real roots if $u \leq 1$, and no positive roots if $1 <u \leq 3$. If $u=3$ then there is only one zero (at the origin), so this case can be discarded.

We consider three cases for the zeros: (i) $\hat{z} = \hat{x}_{m}^{+}>0$ with $u > 3$, (ii) $\hat{z} = -\hat{x}_{m}^{-} \leq 0$ with $u>1$, and (iii) $\hat{z} = -\hat{w}_{m} \in \mathbb{C}$ where $\Re(\hat{w}_{m}) > 0$ and $\Im(\hat{w}_{m}) > 0$. If we include $\hat{z} = -\overline{\hat{w}_{m}}$ this covers all possible zeros, and in all cases the index $m$ is a nonnegative integer, which is bounded in the real cases, and unbounded in the complex case.

We remark that if $u=2n+1$ ($n=1,2,3,\ldots$) we have, to within a positive multiplicative factor function, Hermite polynomials (see (\ref{eq04a})) with only real zeros, and the negative ones following from symmetry from the positive ones. So only the (finite) positive zeros need to be considered, and the formulas presented next in \cref{subsec:poszeros} of course still hold in this important special case.

We find in all three cases the formulas for the uniform asymptotic expansions of the zeros are the same as in \cref{sec:apos}, except that the leading approximations to $z$ and $\zeta$ are changed (but again involve the zeros of certain Airy functions).

\subsection{(i) Positive zeros}
\label{subsec:poszeros}

Here $u>3$, and for each $n=1,2,3,\ldots$ we have from \cite[Sect. 12.11(i)]{NIST:DLMF} that there are $n$ positive roots of (\ref{eq50a}) when $4n-1<u<4n+3$. In the polynomial case $u=2n+1$ there are $n$ zeros, which are symmetric about the origin, with $n/2$ positive ones when $n$ is even and $(n-1)/2$ positive ones when $n$ is odd (and one at the origin in this case). In all cases let $M^{+}=M^{+}(u)$ be the number of positive zeros, which is well-defined from the preceding discussion.

Recall the positive roots of (\ref{eq50a}) are denoted by $\hat{z}=\hat{x}_{m}^{+} > 0$. Then we have from (\ref{eq1.8a}) that these satisfy, in place of (\ref{eq22}),
%%%%%%%%%%
\begin{equation}
\label{eq22a}
 u^{2/3}\mathcal{Z}(u,\hat{x}_{m}^{+})
=\mathrm{a}_{m}
\quad (m=1,2,3,\ldots, M^{+}).
\end{equation}
%%%%%%%%%%

Our desired expansion is again of the form. 
%%%%%%%%%%%%%%%%%%%
\begin{equation}
\label{eq50b}
\hat{x}^{+}_{m} \sim \hat{x}^{+}_{m,0}+\sum_{s=1}^{\infty}
\frac{\hat{x}^{+}_{m,s}}{u^{2s}}
\quad (u \to \infty).
\end{equation}
%%%%%%%%%%%%%%%%%%%
As mentioned above, we only need to find the leading term $\hat{x}^{+}_{m,0}$, since the others then follow from the same formulas as the positive parameter case previously. 

Now for large $u$ we have from (\ref{eq1.11}) and (\ref{eq22a}) $\mathcal{Z}(u,\hat{x}_{m,0}^{+}) \sim \zeta_{m,0}^{+}$ where 
%%%%%%%%%%
\begin{equation}
\label{eq30aa}
\zeta_{m,0}^{+}=\zeta(\hat{x}^{+}_{m,0})
= u^{-2/3}\mathrm{a}_{m} <0
\quad  (m=1,2,3,\ldots ,M^{+}).
\end{equation}
%%%%%%%%%%
Since $\zeta_{m,0}^{+}<0$ it means $\hat{x}_{m,0}^{+} < 1$ (see (\ref{eqZeta1})). Thus in place of $\hat{z}_{m,0}$ as given by (\ref{eq34}) we instead have $\hat{x}_{m,0}^{+} \in (0,1)$ given implicitly by
%%%%%%%%%%
\begin{multline} 
\label{eq32aa}
\arccos\left(\hat{x}_{m,0}^{+}\right)
-\hat{x}_{m,0}^{+}\left(1
-\{\hat{x}_{m,0}^{+}\}^{2}\right)^{1/2}
\\
=\frac{4}{3}(-\zeta_{m,0}^{+})^{3/2}
=\frac{4}{3 u}
\left|\mathrm{a}_{m}\right|^{3/2}
\quad (m=1,2,3,\ldots ,M^{+}).
\end{multline}
%%%%%%%%%%
Note that the LHS is decreasing for $\hat{x}_{m,0}^{+} \in (-1,1)$ and the RHS is increasing as $m$ increases, and as such $\hat{x}_{1,0}^{+}>\hat{x}_{2,0}^{+}>\cdots >\hat{x}_{M^{+},0}^{+}>0$, and the same is true for $\hat{x}_{m}^{+}$. \emph{Thus our zeros in this case are enumerated by decreasing values}.

Note that in (\ref{eqSigma0}) one must replace $\hat{z}_{m,0}$ and $\zeta_{m,0}$ by $\hat{x}^{+}_{m,0}$ and $\zeta_{m,0}^{+}$ respectively, and likewise in (\ref{eqhatz1}), (\ref{eqhatz2}) and the formulas for the subsequent coefficients, which are now labelled $\hat{x}^{+}_{m,s}$ instead of $\hat{z}_{m,s}$ ($s=1,2,3,\ldots$). 

In summary, there are a finite number of positive zeros $\hat{x}_{m}^{+}$ ($m=1,2,3,\ldots, M^{+}$) where $\hat{x}_{m}^{+} > \hat{x}_{m-1}^{+}$, and asymptotically they satisfy (\ref{eq50b}).

\subsection{(ii) Non-positive zeros}

Next consider (\ref{eq50a}) for the non-positive roots $\hat{z}=-\hat{x}_{m}^{-}$ ($\hat{x}_{m}^{-} \geq 0$). We can assume that $u>1$ and $u \neq 2n+1$ ($n=0,1,2,\dots$), as discussed above. Again there are a finite number of these, and denote this number by $M^{-}$. We then enumerate them as follows:
%%%%%%%%%%
\begin{equation}
\label{eq1.8d}
\hat{x}_{1-\vartheta}^{-} > \hat{x}_{2-\vartheta}^{-} > \hat{x}_{3-\vartheta}^{-} \ldots > \hat{x}_{M^{-}-\vartheta}^{-} \geq 0.
\end{equation}
%%%%%%%%%%
Here $\vartheta=\vartheta(u)=1$ if $1<u+2n<\frac{4}{3}$ for some integer $n$, and 0 otherwise, and the reason for introducing this will be become clear shortly.

Unfortunately, unlike $M^{+}$, $M^{-}$ is not generally explicitly given. We do know that $M^{-} \sim M^{+}$ when both are large. We are content to determine $M^{-}=M^{-}(u)$ as the largest positive integer such that our asymptotic expansion for $\hat{x}_{M^{-}-\vartheta}^{-}$ is non-negative, for each given $u$.

Now from (\ref{eq05a}) and (\ref{eq1.8}) (with lower signs in the latter)
%%%%%%%%%%%%%%%%%%%
\begin{equation}
\label{eq1.8b}
U\left(-\tfrac{1}{2}u,-\sqrt {2u} \, \hat{x}_{m}^{-}\right)
=\frac{\Gamma\left(\tfrac{1}{2}u+\tfrac{1}{2}\right)}{\sqrt{2\pi}}\mathcal{Y}(u,\hat{x}_{m}^{-})
\mathcal{A}i\left( u^{2/3}
\mathcal{Z}(u,\hat{x}_{m}^{-}) \right),
\end{equation}
%%%%%%%%%%%%%%%%%%%
where
%%%%%%%%%%
\begin{equation}
\label{eq29b}
\mathcal{A}i(u,z)
=e^{\frac16(3u-1)\pi i}\mathrm{Ai}_{1}(z)
+e^{-\frac16(3u-1)\pi i}\mathrm{Ai}_{-1}(z).
\end{equation}
%%%%%%%%%%
Clearly $\mathcal{A}i(u+2n,z)=(-1)^{n}\mathcal{A}i(u,z)$ ($n \in \mathbb{Z}$), and so all the zeros of this function are determined for $0 \leq u < 2$. We also see that $\mathcal{A}i(u,z)$ is real for $u,z \in \mathbb{R}$. We note from \cite[Eq. 9.2.11]{NIST:DLMF} that when the zeros of $\mathcal{A}i(u,z)$ are the real they are solutions of
%%%%%%%%%%
\begin{equation}
\label{eq29d}
\Re\left\{e^{\frac16(3u-1)\pi i}\mathrm{Ai}_{1}(z)\right\}
=\sin\left(\tfrac{1}{2}u \pi\right) \mathrm{Ai}(z)
+\cos\left(\tfrac{1}{2}u \pi\right) \mathrm{Bi}(z)=0
\quad (z \in \mathbb{R}).
\end{equation}
%%%%%%%%%%

\begin{lemma}
\label{eqlemma2.1}
For $u \in \mathbb{R}$ $\mathcal{A}i(u,z)$ has an infinite number of negative zeros. Moreover, it has one, and only one, positive zero if $1<u+2n<\frac{4}{3}$ for some integer $n$, and no positive zeros otherwise. Finally, it has a zero at $z=0$ iff $u+2n=\frac{4}{3}$ ($n \in \mathbb{Z}$).
\end{lemma}

For $m=1,2,3,\ldots$ let $\mathpzc{a}_{m}^{-}$ denote the $m$th non-positive zero of $\mathcal{A}i(u,z)$, enumerated in the usual manner by increasing absolute values. We label $\mathpzc{a}_{0}$ as the sole positive zero, if the function has one (namely when $u+2n \in (1,\frac43)$, $n=0,1,2,\ldots$). Details on the asymptotics of these zeros is given in \cref{secA}.

Now from (\ref{eq1.8b}) we solve asymptotically for the roots $\hat{x}_{m}^{-}$ of the equation
%%%%%%%%%%%%%%%%%%%
\begin{equation}
\label{eq1.8c}
u^{2/3}\mathcal{Z}(u,\hat{x}_{m}^{-})
=\mathpzc{a}_{m}^{-}
\quad (m =1-\vartheta,2-\vartheta,3-\vartheta,\ldots , M^{-}-\vartheta),
\end{equation}
%%%%%%%%%%%%%%%%%%%
where $\vartheta$ is defined after (\ref{eq1.8d}). When $m=0$ we replace $\mathpzc{a}_{m}^{-}$ by $\mathpzc{a}_{0}$.

The asymptotic expansions are again of the form
%%%%%%%%%%%%%%%%%%%
\begin{equation}
\label{eq50c}
\hat{x}^{-}_{m} \sim \hat{x}^{-}_{m,0}+\sum_{s=1}^{\infty}
\frac{\hat{x}^{-}_{m,s}}{u^{2s}}
\quad (u \to \infty).
\end{equation}
%%%%%%%%%%%%%%%%%%%
The coefficients $\hat{x}^{-}_{m,s}$ ($s=1,2,3,\ldots$) in (\ref{eq50c}) are determined the same way as above, except in place of (\ref{eq34}) we have that the leading terms satisfy
%%%%%%%%%%
\begin{equation}
\label{eq29c}
\zeta_{m,0}^{-}=\zeta(\hat{x}_{m,0}^{-})
= u^{-2/3}\mathpzc{a}_{m}^{-}
\quad  (m=1-\vartheta,2-\vartheta,3-\vartheta,
\ldots ,M^{-}-\vartheta),
\end{equation}
%%%%%%%%%%
again with $\mathpzc{a}_{m}^{-}$ replaced by $\mathpzc{a}_{0}$ if $m=0$.

For $m=1,2,3,\ldots ,M^{-}-\vartheta$ we thus have $\zeta^{-}_{m,0} \leq 0$, therefore $\hat{x}_{m,0}^{-} \leq 1$ and are given by (\ref{eqZeta1}) with $\zeta=\zeta^{-}_{0,0}$ and $\hat{z}=\hat{x}_{m,0}^{-}$. Hence they satisfy
%%%%%%%%%%%%%%%%%%%
\begin{equation} 
\label{eq36}
\arccos\left(\hat{x}_{m,0}^{-}\right)
-\hat{x}_{m,0}^{-}\left(1
-\{\hat{x}_{m,0}^{-}\}^{2}\right)^{1/2}
=\frac{4}{3}(-\zeta_{m,0}^{-})^{3/2}
=\frac{4}{3 u}
\left|\mathpzc{a}_{m}^{-}\right|^{3/2}.
\end{equation}
%%%%%%%%%%%%%%%%%%%

For $1<u+2n < \frac{4}{3}$ ($n \in \mathbb{Z}$) we must also consider $m=0$, and $\hat{x}_{0,0}$ and $\zeta_{0,0}$ correspond to the sole positive zero of $\mathcal{A}i(u,z)$. In this case we have $\zeta^{-}_{0,0} > 0$ from (\ref{eq29c}), and hence $\hat{x}_{0,0}>1$, and so from (\ref{eqZeta}) it instead satisfies
%%%%%%%
\begin{equation}
%\label{eq}
\hat{x}_{0,0}\left(\hat{x}_{0,0}^{2}
-1\right)^{1/2}
-\ln\left\{\hat{x}_{0,0}
+\left(\hat{x}_{0,0}^{2}-1\right)^{1/2}
\right\}
=\frac{4}{3}\zeta_{0,0}^{3/2}
=\frac{4\mathpzc{a}_{0}^{3/2}}{3u}.
\end{equation}
%%%%%%%

\subsection{(iii) Non-real zeros}

Now consider the complex roots of (\ref{eq50a}), which we write as $\hat{z} = -\hat{w}_{m} \in \mathbb{C}$ with $\Re(\hat{w}_{m}) > 0$ and $\Im(\hat{w}_{m}) > 0$. These exist provided $u \neq 2n+1$ for $n=0,1,2,\ldots$, i.e. the non-Hermite polynomial case, which we assume here. By the Schwarz reflection principle $U(-\tfrac{1}{2}u,- \sqrt {2u}\hat{z})$ also has complex zeros in the second quadrant at $\hat{z}=-\overline{\hat{w}_{m}}$.

To examine these we again use (\ref{eq1.8b}), but with $\hat{x}^{-}_{m}$ replaced by $\hat{w}_{m}$. This time we need the corresponding complex zeros of $\mathcal{A}i(u,z)$ in the first and fourth quadrants, and the suitable representation for the former is
%%%%%%%%%%
\begin{equation}
\label{eq52}
\mathcal{A}i(u,z)
=2e^{-\pi i/6} \cos\left(\tfrac{1}{2}u \pi\right)\mathrm{Ai}_{1}(z)
+ie^{-u\pi i/2}\mathrm{Ai}(z),
\end{equation}
%%%%%%%%%%
which comes from (\ref{eqAiConnect}) and (\ref{eq29b}). Note that $\mathrm{Ai}(z)$ is recessive in $\arg(z) < \frac13 \pi$ and $\mathrm{Ai}_{1}(z)$ is recessive in $\frac13 \pi < \arg(z) < \pi$, and hence are a numerically satisfactory pair in the first quadrant.

Denote the non-real zeros of $\mathcal{A}i(u,z)$ in the first quadrant by $z=\mathpzc{a}_{m}$, where again $|\mathpzc{a}_{m}|$ are enumerated so that they are increasing for $m=1,2,3,\ldots$. Observe that by Schwarz symmetry there are also zeros in the fourth quadrant at $z=\overline{\mathpzc{a}_{m}}$. A discussion on the asymptotics of these zeros is given in \cref{secA}.

Corresponding to each of $\mathpzc{a}_{m}$ we have the infinite number of complex roots of (\ref{eq50a}) $\hat{z}=-\hat{w}_{m}$ ($m=1,2,3,\ldots$). Now let $\varsigma_{m}:=\zeta(\hat{w}_{m})$. The equation we again need to solve asymptotically is (\ref{eq1.8c}), but with the real values $\hat{x}_{m}^{-}$ and $\mathpzc{a}_{m}^{-}$ replaced by the complex ones $\hat{w}_{m}$ and $\mathpzc{a}_{m}$, respectively. For $m=1,2,3,\ldots$ the leading approximations to $\hat{w}_{m}$ and $\varsigma_{m}$ are denoted by $\hat{w}_{m,0}$ and $\varsigma_{m,0}$ respectively. Here $\varsigma_{m,0}=u^{-2/3}\mathpzc{a}_{m}$, and $\hat{w}_{m,0}$ is then given by
%%%%%%%%%%
\begin{multline}
\label{eq37}
\hat{w}_{m,0}\left(\hat{w}_{m,0}^{2}-1\right)^{1/2}
-\ln\left\{\hat{w}_{m,0}
+\left(\hat{w}_{m,0}^{2}-1\right)^{1/2}
\right\}
\\
=\frac{4}{3}\varsigma_{m,0}^{3/2}
=\frac{4}{3u}\mathpzc{a}_{m}^{3/2}
\quad (m=1,2,3,\ldots),
\end{multline}
%%%%%%%%%%
where the branches of the logarithm and square roots are such that $\hat{w}_{m,0}$ lies in the first quadrant. In \cref{secA} we show that $|\mathpzc{a}_{m}| \to \infty$ and $\arg(\mathpzc{a}_{m}) \to \frac{1}{3} \pi$ as $m \to \infty$, and hence from (\ref{eq37}) $|\hat{w}_{m,0}| \to \infty$ and $\arg(\hat{w}_{m,0}) \to \frac{1}{4} \pi$ as $m \to \infty$. Numerically it can be more stable to solve this equation by expressing it in the form
%%%%%%%%%%
\begin{multline}
\label{eq37a}
\hat{w}_{m,0}\left[\left(1-
\frac{1}{\hat{w}_{m,0}^{2}}\right)^{1/2}
-\frac{\ln\left(\hat{w}_{m,0}\right)}
{\hat{w}_{m,0}^{2}}
-\frac{1}{\hat{w}_{m,0}^{2}}
\ln\left\{1
+\left(1-\frac{1}{\hat{w}_{m,0}^{2}}\right)^{1/2}
\right\} \right]^{1/2}
\\
=\frac{2}{\sqrt{3u}}\mathpzc{a}_{m}^{3/4}
\quad (m=1,2,3,\ldots).
\end{multline}
%%%%%%%%%%

With these initial coefficients we proceed as in all the previous cases to determine the subsequent coefficients in the expansion, $\hat{w}_{m,s}$ say, and our desired asymptotic expansion is then given by
%%%%%%%%%%%%%%%%%%%
\begin{equation}
\label{eq50d}
\hat{w}_{m} \sim \hat{w}_{m,0}+\sum_{s=1}^{\infty}
\frac{\hat{w}_{m,s}}{u^{2s}}
\quad (u \to \infty).
\end{equation}
%%%%%%%%%%%%%%%%%%%

In order to check the accuracy of our zero approximations we could compute the left hand side of the following expression (see \cite[Eq. 12.2.18]{NIST:DLMF}) using our approximation for $\hat{w}_{m}$
%%%%%%%%%%%%%%%%%%%
\begin{equation}
\label{eq67}
1+\frac{ie^{-u\pi i/2}U\left(\tfrac{1}{2}u,  
i\sqrt {2u} \,\hat{w}_{m}\right)}
{U\left(\tfrac{1}{2}u, -i \sqrt {2u} 
\, \hat{w}_{m}\right)}=0.
\end{equation}
%%%%%%%%%%%%%%%%%%%
Additional tests are described in the next section for real and complex zeros.

\section{Numerical results}
\label{num:tests}

Here we numerically test the accuracy of our uniform expansions in both the real variable case ($a$ negative corresponding to Hermite polynomials) and complex variable case ($a$ both positive and negative).

\subsection{Real zeros: Hermite polynomial case}

We consider the expansion (\ref{eq50b}) for the case $a=-\frac{61}{2}$ ($u=61$ in (\ref{eq50a})). Thus from (\ref{eq04a}) we are considering the positive zeros $x=\sqrt{61}\,\hat{x}^{+}_{m}$ ($m=1,2,3,\ldots$ 15) of the Hermite polynomial $H_{30}(x)$.

Consider
%%%%%%%%%%%%%%%%%%%
\begin{equation}
\label{eq4.1}
V_{n}(x)=n! \,2^{\frac12 n}
e^{\frac12 x^2}
V\left(-n - \tfrac12, 2^{\frac12}x\right),
\end{equation}
%%%%%%%%%%%%%%%%%%%
where $V(a,z)$ is a companion solution to $U(a,z)$ of Weber's differential equation (\ref{eq01}) (see \cite[sec. 12.5]{NIST:DLMF} for various integral representations). For real $x$ $V_{n}(x)$ represents a numerically satisfactory companion to $H_{n}(x)$ that satisfies the same differential equation, namely
%%%%%%%%%%%%%%%%%%%
\begin{equation}
\label{eq4.2}
\frac{d^2 y}{dx^2}-2x\frac{dy}{dx}+2ny=0.
\end{equation}
%%%%%%%%%%%%%%%%%%%

\begin{figure}[h!]
 \centering
 \includegraphics[
 width=0.7\textwidth,keepaspectratio]{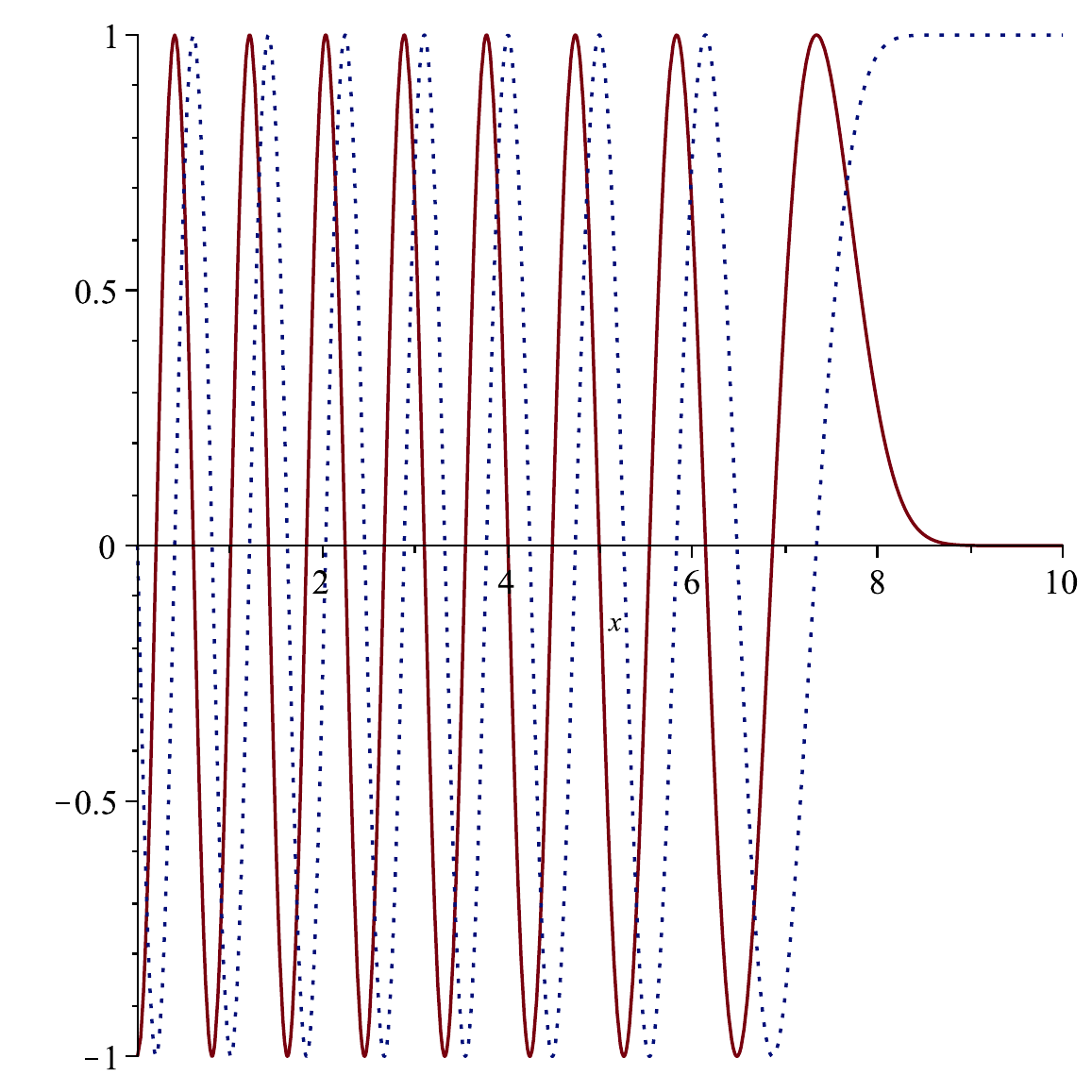}
 \caption{Graphs of $H_{30}(x)/F_{30}(x)$ (solid curve) and $V_{30}(x)/F_{30}(x)$ (dotted curve)}
 \label{fig:Hermite}
\end{figure}

To illustrate this, consider the envelope (c.f. \cite[Eqs. 12.2.21	and 14.15.23]{NIST:DLMF})
%%%%%%%%%%%%%%%%%%%
\begin{equation}
\label{eq4.3}
F_{n}(x)=\left\{H_{n}^{2}(x)+V_{n}^{2}(x)\right\}^{1/2}.
\end{equation}
%%%%%%%%%%%%%%%%%%%
Then the graphs of $H_{30}(x)/F_{30}(x)$ (solid curve) and $V_{30}(x)/F_{30}(x)$ (dotted curve) are depicted in \cref{fig:Hermite} for $0 \leq x \leq 10$. We observe that both (scaled) functions have amplitudes of oscillation approximately equal to unity, and are approximately $\frac12 \pi$ out of phase in the interval where they oscillate. Due to the amplitude of oscillation being close to one the value of $H_{30}(x)/F_{30}(x)$ at a truncated asymptotic expansion for a zero gives a good estimate on the accuracy of said approximation, as we examine next.

Let us label the positive zeros of $H_{n}(x)$ by $x=h_{n,m}$ for $m=1,2,3,\ldots \lfloor{n/2}\rfloor$. Recall from \cref{subsec:poszeros} that these zeros are enumerated by decreasing values. Then from (\ref{eq04a}), (\ref{eq50a}) and (\ref{eq50b}), and taking five terms in the expansion, we have the approximation
%%%%%%%%%%%%%%%%%%%
\begin{equation}
%\label{eq50b}
h_{30,m} \approx \varrho_{m}
:= \sqrt{61} \, \sum_{s=0}^{4}
\frac{\hat{x}^{+}_{m,s}}{61^{2s}}
\quad (m=1,2,3,\ldots 15).
\end{equation}
%%%%%%%%%%%%%%%%%%%
From this we find, for example, $h_{30,1}=6.86334\cdots$ and $h_{30,15}=0.20112\cdots$. 

\begin{table}[h!]
\centering
$\displaystyle
\begin{array}{cc|cc|cc}
\hline
   m & |H_{30}(\varrho_{m})/F_{30}(\varrho_{m})| & m 
&  |H_{30}(\varrho_{m})/F_{30}(\varrho_{m})|
& m &  |H_{30}(\varrho_{m})/F_{30}(\varrho_{m})|
\\ 
  \hline
1 &  2.075\times 10^{-17}  &  6 &  1.439 \times 10^{-16} 
& 11  & 1.006 \times 10^{-15}  \\
2 &  3.465 \times 10^{-17}  & 7 &  2.051 \times 10^{-16} 
& 12  & 1.584 \times 10^{-15} \\
3 &  5.103 \times 10^{-17}  & 8 &  2.965 \times 10^{-16}
& 13  & 2.563 \times 10^{-15} \\
4 &  7.248 \times 10^{-17}  & 9  & 4.362 \times 10^{-16} 
& 14  &  4.275 \times 10^{-15} \\
5 &  1.020 \times 10^{-16}  & 10  &  6.550 \times 10^{-16}  
& 15  &  7.372 \times 10^{-15} \\
     \end{array}
$
\caption{Values of $|H_{30}(\varrho_{m})/F_{30}(\varrho_{m})|$ for $m=1,2,3,\ldots 15$.}
\label{table1}
\end{table}

In \cref{table1} we evaluate $|H_{30}(\varrho_{m})/F_{30}(\varrho_{m})|$ for $m=1,2,3,\ldots 15$, which illustrates the uniform accuracy of all our asymptotic expansions. Note that the accuracy improves as the size of the zero increases (correspondingly as $m$ decreases), which demonstrates the doubly uniform nature of our approximations, namely being valid as the magnitude of the zeros approaching infinity as well as $|a| \to \infty$.

\subsection{Complex zeros}

In order to avoid wrong branches when computing complex values of $\zeta$, we recast (\ref{eqZeta}) and (\ref{eqZeta1}) into forms that are more appropriate, namely
%%%%%%%%%%%%%%%%%%%
\begin{equation}
\label{eq23}
\zeta=\hat{z}^{4/3} \left[\tfrac{3}{4}\left\{
\left(1-\hat{z}^{-2}\right)^{1/2}
-\hat{z}^{-2}\ln \left\{ 1+
\left(1-\hat{z}^{-2}\right)^{1/2}
\right\} 
-\hat{z}^{-2}\ln (\hat{z})\right\}\right]^{2/3},
\end{equation}
%%%%%%%%%%%%%%%%%%%
for $|\hat{z}| \geq 1$, and 
%%%%%%%%%%%%%%%%%%%
\begin{equation}
\label{eq24}
\zeta=-\left[\tfrac{3}{4}\left\{
\arccos(\hat{z})-\hat{z}
\left(1-\hat{z}\right)^{1/2}
\right\}\right]^{2/3},
\end{equation}
%%%%%%%%%%%%%%%%%%%
for $|\hat{z}| < 1$. In both (\ref{eq23}) and (\ref{eq24}) principal branches are taken for all multi-valued terms, and in particular for the inverse cosine we have
%%%%%%%%%%%%%%%%%%%
\begin{equation}
\label{eq25}
\arccos(\hat{z}) =-i\ln  \left\{\hat{z}
+i\left(1-\hat{z}\right)^{1/2} \right\}.
\end{equation}
%%%%%%%%%%%%%%%%%%%

For $|\hat{z}| \geq 1$ ($\hat{z} \neq 1$) we likewise compute $\beta$ by recasting (\ref{eq06}) into the form
%%%%%%%%%%%%%%%
\begin{equation}
\label{eq26}
\beta=\left(1-\hat{z}^{-2}\right)^{-1/2},
\end{equation}
%%%%%%%%%%%%%%%
and for $|\hat{z}|<1$ with $\Im(\hat{z}) \gtrless 0$
%%%%%%%%%%%%%
\begin{equation} 
\label{eq27}
\beta=\mp i \hat{z} \left(1-\hat{z}^{-2}\right)^{-1/2}.
\end{equation}
%%%%%%%%%%%%%%%%%%%

For testing the accuracy of the uniform asymptotic expasions for the zeros of $U(a,z)$, we use the results given in \cite{Segura:2013:CCZ}. 
In that reference, a method for finding the complex zeros of solutions of second order ordinary differential equations
is described. The method makes use first of a qualitative analysis of the approximate Liouville-Green Stokes lines (SLs) and
anti-Stokes lines (ASLs). As discussed in  \cite{Segura:2013:CCZ}, the structure of the exact zeros will follow very closely
the ASLs. This qualitative analysis is combined with the application of a  fixed point method $z_{n+1}=T(z_n)$
(which has fourth-order convergence and good non-local behaviour) and the use of carefully
selected step functions (displacements) $H^{\pm}(z)$.

To compute the complex zeros of $U(a,z)$, the following displacements  $H^{+}(z)$ and iterating function $T(z)$ 
are used:
%%%%%%%%%%%%%%%%%%%
\begin{equation}
H^+(z)=z+\pi  \left(-\tfrac14z^{2}-a\right)^{-1/2},
\label{displa}
\end{equation}
%%%%%%%%%%%%%%%%%%%
and
%%%%%%%%%%%%%%%%%%%
\begin{equation}
T(z)=z-\left(-\tfrac14z^{2}-a\right)^{-1/2}
\arctan\left\{
\left(-\tfrac14z^{2}-a\right)^{1/2}
\frac{U(a,z)}{U'(a,z)}
\right\}.
\label{iter}
\end{equation}
%%%%%%%%%%%%%%%%%%%
In the comparisons, we test the accuracy of the approximations to the zeros $z=z_m$ of $U(a,z)$ using the following functions

\begin{equation}
\begin{array}{l}
g_1=\left[\{\Re (z_m)\}^2 + \{\Im (z_m)\}^2 \right]^{1/2} \,,\\
\\
g_2= \Re (z_m)/\Im (z_m) \quad (\Re (z_m)\Im (z_m) \neq 0)\,.
\end{array}
\end{equation}
In the notation of \cref{sec:apos,sec:aneg} we observe from (\ref{eq21}) and (\ref{eq50a}) that $z_{m}=2 i \sqrt{a}\,\hat{z}_{m}$ for $a>0$, and $z_{m}=\pm 2 \sqrt{|a|}\,\hat{x}_{m}^{\pm}$ or $z_{m}=-2 \sqrt{|a|}\,\hat{w}_{m}$ for $a<0$ (according to whether $z_{m}$ is positive, negative, or complex).

Then, we evaluate the relative errors

 \begin{equation}
 \epsilon_1=\left|1-\displaystyle\frac{g_{1}^{\mathrm{approx} }}{g_{1}^{\mathrm{numer}}} \right|,\,\,
\epsilon_2=\left|1-\displaystyle\frac{g_{2}^{\mathrm{approx}}}{g_{2}^{\mathrm{numer}}} \right|,
\end{equation}
where $g_{j}^{\mathrm{approx}}$ ($j=1,2$) are evaluated using the approximations given in previous sections and  
$g_{j}^{\mathrm{numer}}$ ($j=1,2$) are computed with our numerical method implemented in Maple with a large number of digits. For a fixed precision implementation of the algorithm, the computational scheme suggested in \cite{Dunster:2024:CPC} could be used to evaluate the function $U(a,z)$ needed in the iterating function $T(z)$ \eqref{iter}
(the derivative $U'(a,z)$ can be computed using, for example, the recurrence relation \cite[Eq. 12.8.2]{NIST:DLMF}).

\begin{table}[h!]
\begin{adjustbox}{width=\textwidth}
$\displaystyle
\begin{array}{cccc}
\hline
   m & z_m & \epsilon_1 &  \epsilon_2 \\ 
  \hline
           1 &-{\bf 1.382736145125}9055+{\bf 6.60363420332}86323i         &   2.3 \times 10^{-13}    & 3.9  \times 10^{-13} \\
           2 &   -{\bf 2.36697098755}73483+{\bf 7.25076501051}86024i       &    1.3 \times 10^{-13}      & 2.4   \times 10^{-13}   \\
           3 &   -{\bf 3.14303439319}50775+{\bf 7.7865053482}195365i      &   8.4 \times 10^{-14}       & 1.8  \times 10^{-13}    \\
           4 &   -{\bf 3.808424713323}3240+{\bf 8.262102283248}3978i      &  6.0 \times 10^{-14}        & 1.4  \times 10^{-13}        \\
           5 &   -{\bf 4.40113226187}31031+{\bf 8.697352864671}4638i      &   4.6 \times 10^{-14}   &    1.1 \times 10^{-13}      \\
    50 &   -{\bf 16.825271666405}126 +{\bf 19.292382093177}420i      &  1.4 \times 10^{-15}        & 1.3  \times 10^{-14}        \\
     100 &   -{\bf 24.3108724465970}90  +{\bf  26.292345765760}354i      &   4.4\times 10^{-16}   &    6.8 \times 10^{-15}      \\
     \end{array}
$
\end{adjustbox}
\caption{\label{table2} Test for the asymptotic expansions of the zeros of $U(8.3,z)$.}
\end{table}

\begin{table}[h!]
\begin{adjustbox}{width=\textwidth}
$\displaystyle
\begin{array}{cccc}
\hline
   m & z_m & \epsilon_1 &  \epsilon_2 \\ 
  \hline
           1 & -{\bf 1.206751169454753}4+{\bf 9.729142195621040}3i       & 4.4 \times    10^{-17}    &  6.9 \times      10^{-17}     \\
           2 &   -{\bf 2.085091237010430}7+{\bf 10.277292389190367}i    &   3.2 \times 10^{-17}   &  5.3 \times 10^{-17}  \\
           3 & -{\bf 2.788861620217136}1+{\bf 10.731269264892200}i    &   2.5 \times  10^{-17}       & 4.3 \times   10^{-17}    \\
           4 &   -{\bf 3.399747162700204}1+{\bf 11.135489161113063}i     &  2.0 \times   10^{-17}       & 3.6 \times  10^{-17}    \\
           5 &   -{\bf 3.949364339009171}2+{\bf 11.506895318690518}i     & 1.6 \times    10^{-17}   &  3.0 \times  10^{-17}   \\
50 &   -{\bf   16.118357080255495  }+{\bf 21.073613351807242   }i     & 9.3\times    10^{-18}   &  4.1 \times  10^{-18}   \\
100 &   -{\bf  23.642327373211272 }+{\bf 27.734831831550747  }i     & 2.9 \times    10^{-19}   &  2.1 \times  10^{-18}   \\
     \end{array}
$
\end{adjustbox}
\caption{\label{table3} Test for the asymptotic expansions of the zeros of $U(20.3,z)$.
}
\end{table}

\begin{table}[h!]
\begin{adjustbox}{width=\textwidth}
$\displaystyle
\begin{array}{cccc}
\hline
   m & z_m & \epsilon_1 &  \epsilon_2 \\ 
  \hline
           1 &   -{\bf  5.690558573}7972570+{\bf 1.3832406806}543917i           &        3.8 \times 10^{-12}    &   6.7 \times 10^{-12} \\
           2 &   -{\bf  6.42034330496}08671+{\bf 2.4184037014}614955i        &  1.8   \times  10^{-12}        & 3.8 \times 10^{-12}     \\
           3 &      -{\bf 7.00528370942}20902+{\bf 3.2229279813}213036i      &    1.1   \times 10^{-12}      &   2.6 \times  10^{-12}    \\
           4 &      -{\bf 7.51760677349}16861+{\bf 3.90724536328}57412i      &  7.8   \times  10^{-13}        &   1.9 \times 10^{-12}         \\
           5 & -{\bf 7.9826003951}377883+{\bf 4.5135383156}131224i      &  5.8   \times  10^{-13}    &   1.6  \times  10^{-12}       \\
 50 &      -{\bf  19.07513238506}2910       +{\bf  17.16307450067}4282i   &  1.5   \times  10^{-14}        &   1.7 \times 10^{-13}         \\
       100 & -{\bf    25.98902178504}7848  +{\bf   24.45308013876}8002i   &  4.7  \times  10^{-15}    &   9.2 \times  10^{-14}       \\
     \end{array}
$
\end{adjustbox}
\caption{\label{table4} Test for the asymptotic expansions of the zeros of $U(-6.2,z)$.
}
\end{table}

\begin{table}[h!]
\begin{adjustbox}{width=\textwidth}
$\displaystyle
\begin{array}{cccc}
\hline
   m & z_m & \epsilon_1 &  \epsilon_2 \\ 
  \hline
           1 &   -{\bf 5.690558573810462967}2+{\bf 1.3832406806482687014}i        &        8.6 \times 10^{-20}    &   2.9 \times 10^{-20} \\
           2 &    -{\bf 6.420343304969841599}5+{\bf 2.4184037014557299517}i      &  3.5  \times  10^{-20}        & 5.1 \times  10^{-21}     \\
           3 &   -{\bf 7.00528370942923144}89+{\bf 3.222927981316204036}7i    &    5.8   \times  10^{-20}      &   2.9 \times  10^{-20}    \\
           4 &   -{\bf 7.51760677349780159}47+{\bf 3.9072453632811518184}i      &  9.2   \times  10^{-20}        &   3.6 \times  10^{-20}         \\
           5 &   -{\bf 7.982600395143232619}5+{\bf 4.513538315608912947}3i   &  8.7   \times  10^{-20}    &   1.6  \times  10^{-20}       \\
 50 &   -{\bf   19.075132385064583145 }+{\bf 17.16307450067271792   }4i      &  6.7   \times  10^{-21}        &   2.1 \times  10^{-20}         \\
     100 &   -{\bf  25.989021785049034971 }+{\bf   24.453080138766863354}i   &  1.7   \times  10^{-21}    &   1.0  \times  10^{-21}       \\

     \end{array}
$
\end{adjustbox}
\caption{\label{table5} Test for the asymptotic expansions of the zeros of $U(-6.2,z)$  with an improvement in accuracy similar to the one described for the zeros of Airy functions in the appendix.
}
\end{table}

In \cref{table2,table3} we show examples (for $a=8.3,\,20.3$) of the results obtained with the approximations given in \cref{sec:apos}. We use $5$ terms in the expansion  \eqref{eq50}. The second column show the values of the approximations to the zeros for different values of $m$.
The digits in bold are the ones that match those obtained with the numerical computation.
The last two columns in the tables show the relative errors $\epsilon_1$ and $\epsilon_2$ obtained in the comparisons with the numerical
method. 

As observed, the approximations yield very good accuracy even for small values of  $m$: for $m \le 5$, the accuracy is close to $10^{-14}$ for $a=8.3$ and better than $10^{-17}$ for
$a=20.3$. Also, as expected the accuracy of the approximations improves as $m$ increases.

In \cref{table4} we show an example for $a$ negative ($a=-6.2$) and $z_{m}=-\hat{w}_{m} \in \mathbb{C}$. As before, we use $5$ terms, this time in the expansion (\ref{eq50d}). \cref{table5} also shows the results of the asymptotic expansion of the zeros of  $U(-6.2,z)$, but now additionally using an improvement in accuracy of the zeros of Airy functions via the method described at the end of \cref{secA}.

\begin{figure}[H]
 \centering
 \includegraphics[
 width=0.6\textwidth,keepaspectratio]{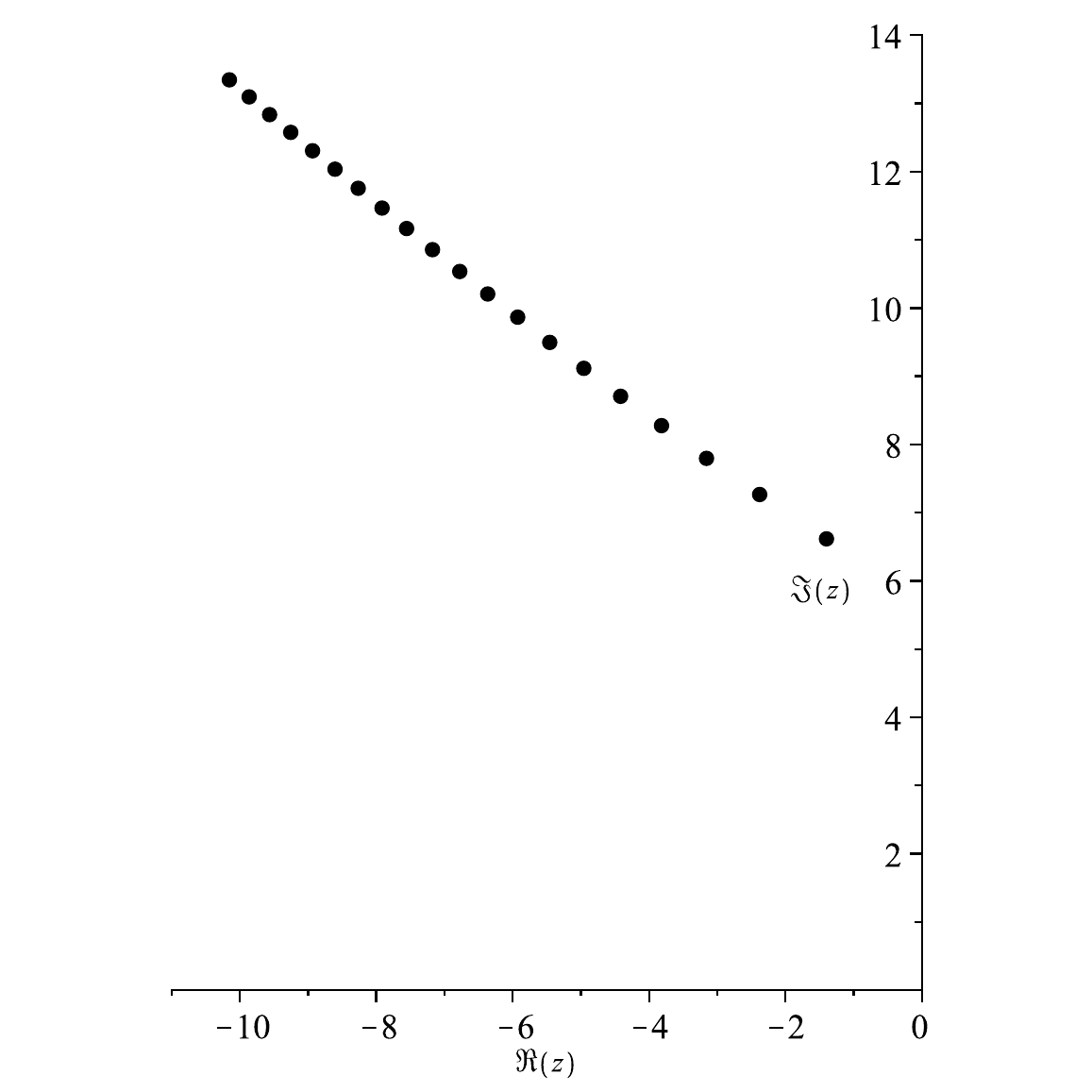}
 \caption{Plot of first $20$ zeros of $U(8.3,z)$ in the second quadrant obtained from the asymptotic approximations.}
 \label{fig:UposaZeroPlot}
\end{figure}

\begin{figure}[H]
 \centering
 \includegraphics[
 width=0.6\textwidth,keepaspectratio]{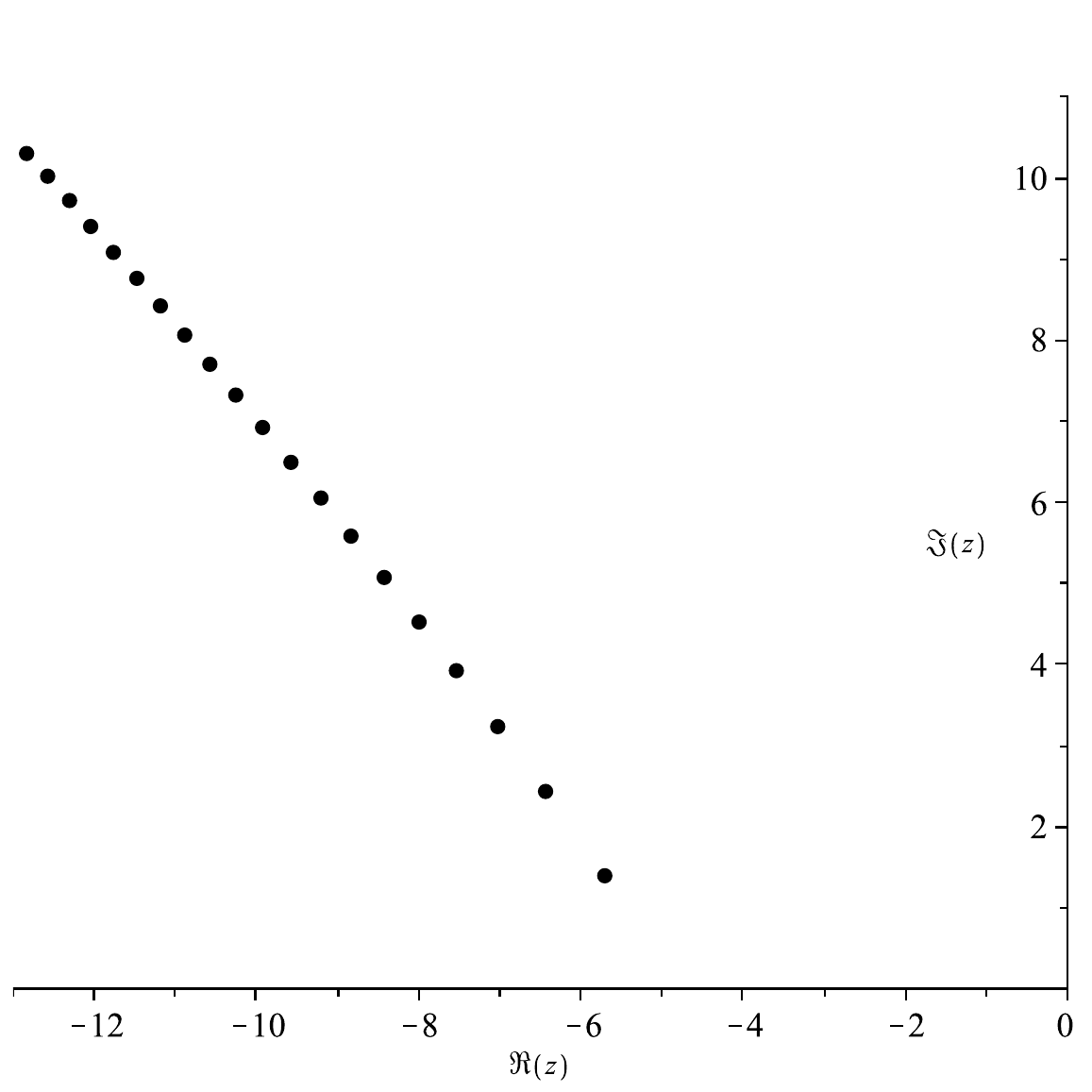}
 \caption{Plot of first $20$ zeros of $U(-6.2,z)$ in the second quadrant obtained from the asymptotic approximations.}
 \label{fig:UnegaZeroPlot}
\end{figure}

\begin{figure}[H]
 \centering
 \includegraphics[
 width=0.6\textwidth,keepaspectratio]{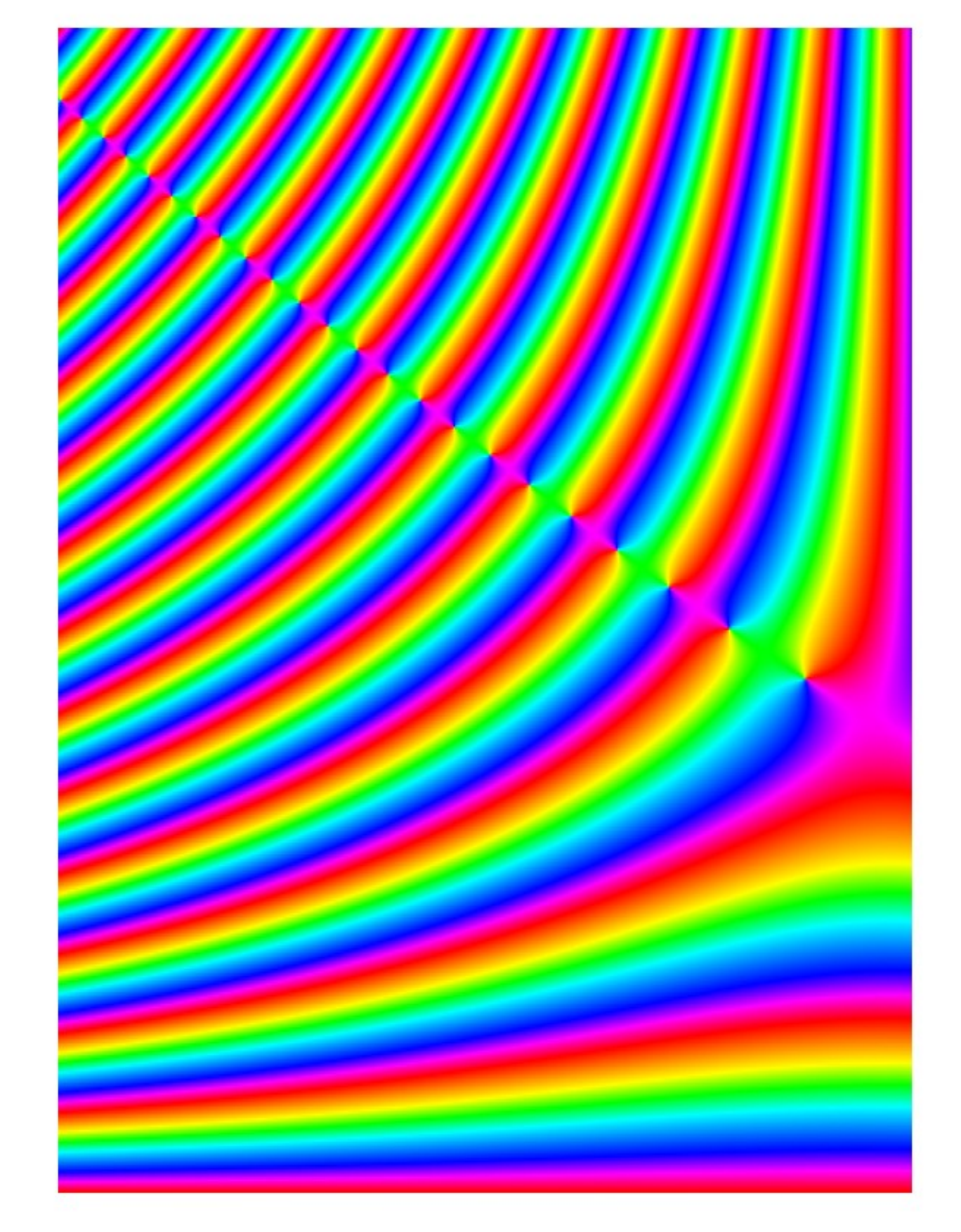}
 \caption{Phase plot of $U(8.3,z)$ in the second quadrant where the first $20$ zeros are located.}
 \label{fig:phaseplotapos}
\end{figure}

\begin{figure}[H]
 \centering
 \includegraphics[
 width=0.6\textwidth,keepaspectratio]{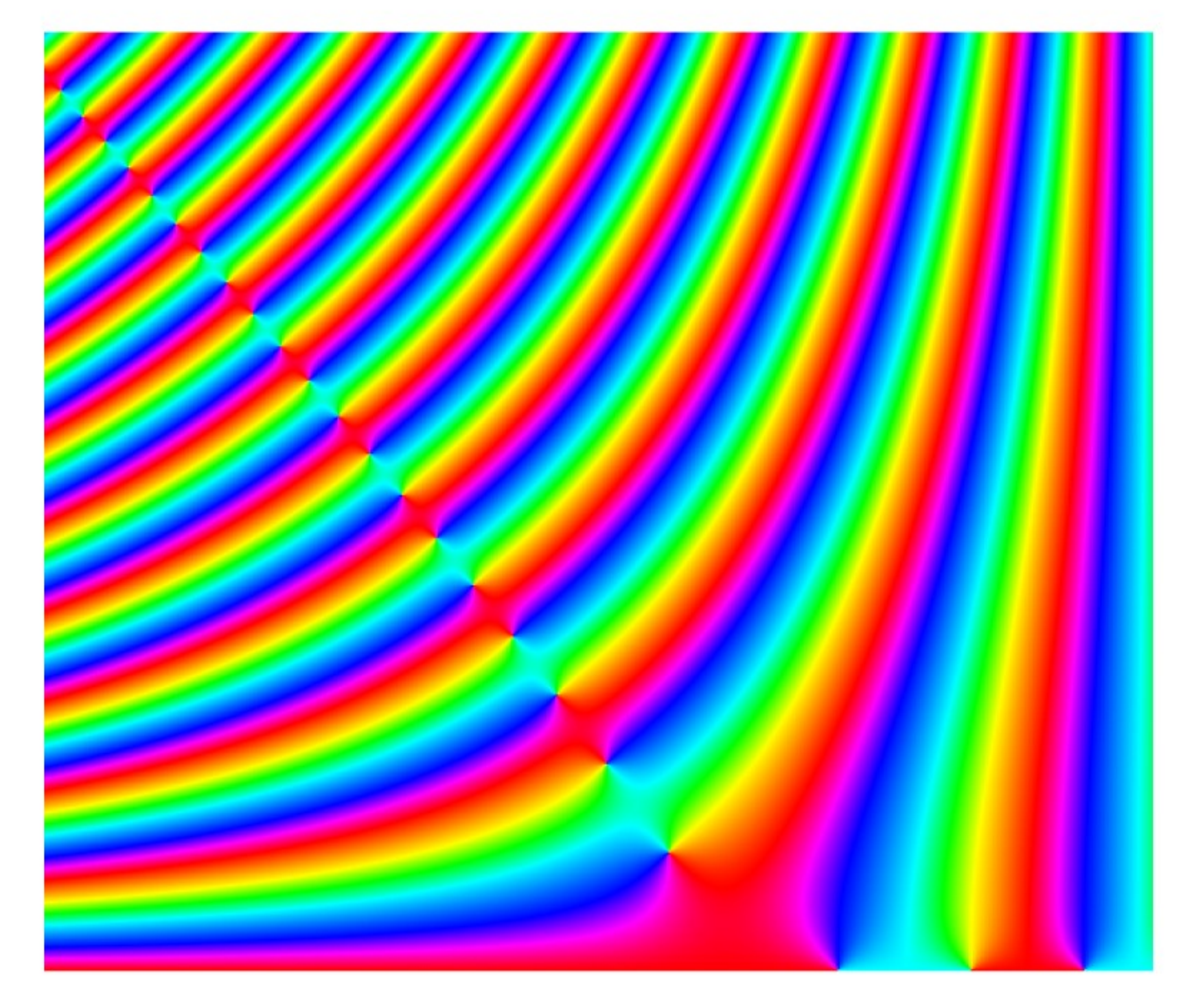}
 \caption{ Phase plot of $U(-6.2,z)$ in the second quadrant where the first $20$ zeros are located.   }
 \label{fig:phaseplotaneg}
\end{figure}

In \cref{fig:UposaZeroPlot,fig:UnegaZeroPlot} the first $20$ zeros of $U(8.3,z)$ and $U(-6.2,z)$, respectively, in the second quadrant are plotted using our asymptotic approximations. \cref{fig:phaseplotapos,fig:phaseplotaneg} show colored representations of the phase of the functions $U(8.3,z)$ and $U(-6.2,z)$, respectively, again in the region of the second quadrant where the first $20$ zeros obtained with the asymptotic approximations are located. As can be seen, the positions of the zeros coincide with the obtained approximations. The values of $U(a,z)$ to obtain the phase representations have been obtained using the methods given in \cite{Dunster:2024:CPC}.

\appendix
\section{Zeros of Airy functions}
\label{secA}

A discussion on the asymptotic behaviour of $\mathpzc{a}_{m}^{-}$ as $m \to \infty$ is given in \cite[Eqs. (21) and (22)]{Gil:2014:CZA}, where in Eq. (11) of that reference we identify $ \alpha=\frac12 (1-u)\pi$. Note that, in terms of their notation, the smallest value of $k$ is not necessarily $1$. In our case $m=1$ is the index for the largest non-positive zero, and with this in mind let us define a periodic function
%%%%%%%%%%%%%%%%%%%
\begin{equation}
\label{eq29e}
\mu(u)=
    \begin{cases}
        2u & (0 \leq u < \tfrac{4}{3})\\
        2u-4 & (\tfrac{4}{3} \leq u \leq 2)
    \end{cases},
\end{equation}
%%%%%%%%%%%%%%%%%%%
with $\mu(u+2n)=\mu(u)$ for $n=1,2,3,\ldots$. In (\ref{eq29e}) the cut off value $u=\frac{4}{3}$ corresponds to the emergence of a first non-positive zero at $\mathpzc{a}_{1}^{-}=0$, necessitating a shift in the index for all the other zeros. Thus we have the asymptotic expansion
%%%%%%%%%%%%%%%%%%%
\begin{equation}
\label{eq29g}
\mathpzc{a}_{m}^{-}=
-T\left(\tfrac{3}{8} \pi \tau_{m}\right)
\quad (m=2,3,4,\ldots),
\end{equation}
%%%%%%%%%%%%%%%%%%%
where
%%%%%%%%%%%%%%%%%%%
\begin{equation}
\label{eq29f}
\tau_{m}=4m-3+\mu(u),
\end{equation}
%%%%%%%%%%%%%%%%%%%
and as $t \to \infty$
%%%%%%%%%%%%%%%%%%%
\begin{equation}
\label{eq62}
T(t)\sim t^{2/3}\left(1+\frac{5}{48 \, t^{2}}
-\frac{5}{36 \, t^{4}}
+\frac{77125}{82944 \, t^{6}}
-\frac{108056875}{6967296 \, t^{8}}
-\cdots\right).
\end{equation}
%%%%%%%%%%%%%%%%%%%

Next consider the asymptotics of $\mathpzc{a}_{m}$. Define $\eta=\frac{2}{3}z^{3/2}$ with the principal root taken. Then from \cite[Thm. 2.4]{Dunster:2021:SEB} we have as $z \to \infty$ for $\delta>0$
%%%%%%%%%%%%%%%%%%%
\begin{equation}
\label{eq54}
\mathrm{Ai}(z) 
\sim \frac{1}{2\pi^{1/2}z^{1/4}}
\exp \left\{ -\eta
+\sum\limits_{s=1}^{\infty}{(-1)^{s}
\frac{a_{s}}{s\eta^{s}}}\right\}
\quad (|\arg(z) |\leq \pi - \delta),
\end{equation}
%%%%%%%%%%%%%%%%%%%
and
%%%%%%%%%%%%%%%%%%%
\begin{equation}
\label{eq55}
\mathrm{Ai}_{\pm 1}(z) 
\sim  \frac{e^{\pm \pi i/6}}{2\pi^{1/2}
z^{1/4}}\exp \left\{ \eta
+\sum\limits_{s=1}^{\infty}\frac{a_{s}}
{s \eta^{s}}\right\}
\quad \left(\left|\arg\left(z e^{\mp 2\pi i/3}\right)
\right| \leq \pi - \delta\right).
\end{equation}
%%%%%%%%%%%%%%%%%%%

Consider first $\cos(\frac{1}{2}u \pi)>0$. From (\ref{eq52}) we wish to asymptotically solve for $z$ in the first quadrant
%%%%%%%%%%%%%%%%%%%
\begin{equation}
\label{eq53}
e^{-\pi i/6} e^{\rho}\mathrm{Ai}_{1}(z)
+e^{-\rho}\mathrm{Ai}(z)=0,
\end{equation}
%%%%%%%%%%%%%%%%%%%
where
%%%%%%%%%%%%%%%%%%%
\begin{equation}
\label{eq53a}
\rho=\tfrac{1}{2}\ln\left\{
2\cos\left(\tfrac{1}{2}u \pi\right)\right\}
+\tfrac{1}{4}(u-1)\pi i.
\end{equation}
%%%%%%%%%%%%%%%%%%%
Thus as $z \to \infty$ for $m=1,2,3,\ldots$ and some integer $m^{+}$ we have from (\ref{eq54}) - (\ref{eq53a})
%%%%%%%%%%%%%%%%%%%
\begin{equation}
\label{eq56}
\eta+\sum\limits_{s=0}^{\infty}
\frac{a_{2s+1}}{(2s+1) \eta^{2s+1}}
\sim \frac{1}{4}\left(4m+4m^{+}-u-1\right)\pi i
-\frac{1}{2}\ln\left\{
2\cos\left(\tfrac{1}{2}u \pi\right)\right\}.
\end{equation}
%%%%%%%%%%%%%%%%%%%

As we shall see, for each fixed $u$ we require for $m=1,2,3,\ldots$ that $4m+4m^{+}-u-1> 0$ in order for the non-real zeros to be in the first quadrant. Thus take $m^{+}$ to be the smallest integer such that $4m^{+}> u-3$, which gives 
%%%%%%%%%%%%%%%%%%%
\begin{equation}
\label{eq57}
m^{+}= \lfloor \tfrac{1}{4}(u+1) \rfloor.
\end{equation}
%%%%%%%%%%%%%%%%%%%

Now compare (\ref{eq56}) with the expansion for the negative zeros of $\mathrm{Ai}(z)$, which from \cite[Eq. 9.2.14]{NIST:DLMF} and (\ref{eq54})
%%%%%%%%%%%%%%%%%%%
\begin{equation}
\label{eq58}
\mathrm{Ai}(-z)
=e^{\pi i/3}\mathrm{Ai}\left(ze^{\pi i/3}\right)
+e^{-\pi i/3}\mathrm{Ai}\left(ze^{-\pi i/3}\right),
\end{equation}
%%%%%%%%%%%%%%%%%%%
gives for $m=1,2,3,\ldots$ 
%%%%%%%%%%%%%%%%%%%
\begin{equation}
\label{eq59}
\eta-\sum\limits_{s=0}^{\infty} (-1)^{s}
\frac{a_{2s+1}}{(2s+1) \eta^{2s+1}}
\sim \frac{1}{4}\left(4m-1\right)\pi,
\end{equation}
%%%%%%%%%%%%%%%%%%%
where again $\eta=\frac{2}{3}z^{3/2}$. Letting $M=4m-1 \to \infty$ and inverting gives an asymptotic expansion for $\eta$, and then using $z=(3\eta/2)^{2/3}$ and expanding this as another asymptotic expansion in inverse powers of $M$ yields the well-known asymptotic expansion for the zeros $\mathrm{a}_{m}$ of $\mathrm{Ai}(x)$ (see for example, \cite[Eqs. 9.9.6 and 9.9.18]{NIST:DLMF}).

If in (\ref{eq59}) we replace $\eta$ by $-i \eta$ (so that $z$ is replaced by $ze^{-\pi i/3}$), and $4m-1$ by $\tau_{m}^{+}$ where
%%%%%%%%%%%%%%%%%%%
\begin{equation}
\label{eq60}
\tau_{m}^{+}
= 4m+4m^{+}-u-1
+\frac{2i}{\pi}\ln\left\{
2\cos\left(\frac{1}{2}u \pi\right)\right\},
\end{equation}
%%%%%%%%%%%%%%%%%%%
we obtain (\ref{eq56}). Thus the asymptotic expansion for $\mathpzc{a}_{m}$ is given by \cite[Eqs. 9.9.6 and 9.9.18]{NIST:DLMF} when modified in a similar way. Specifically, 
%%%%%%%%%%%%%%%%%%%
\begin{equation}
\label{eq61}
\mathpzc{a}_{m}=e^{\pi i/3}T\left(
\tfrac{3}{8} \pi \tau_{m}^{+}\right),
\end{equation}
%%%%%%%%%%%%%%%%%%%
where $T(t)$ is given by (\ref{eq62}). From (\ref{eq62}), (\ref{eq60}) and (\ref{eq61}) it is now evident why we require $4m+4m^{+}-u-1>0$ for zeros in the first quadrant.

Consider now $\cos(\frac{1}{2}u \pi)<0$. Then in place of (\ref{eq56}) we find in a similar manner
%%%%%%%%%%%%%%%%%%%
\begin{equation}
\label{eq63}
\eta+\sum\limits_{s=0}^{\infty}
\frac{a_{2s+1}}{(2s+1) \eta^{2s+1}}
\sim \frac{1}{4}\left(4m+4m^{-}-u+1\right)\pi i
-\frac{1}{2}\ln\left\vert
2\cos\left(\tfrac{1}{2}u \pi\right)\right\vert,
\end{equation}
%%%%%%%%%%%%%%%%%%%
where 
%%%%%%%%%%%%%%%%%%%
\begin{equation}
\label{eq64}
m^{-}= \lfloor \tfrac{1}{4}(u-1) \rfloor.
\end{equation}
%%%%%%%%%%%%%%%%%%%
Then we arrive at
%%%%%%%%%%%%%%%%%%%
\begin{equation}
\label{eq65}
\mathpzc{a}_{m}=e^{\pi i/3}T\left(
\tfrac{3}{8} \pi \tau_{m}^{-}\right),
\end{equation}
%%%%%%%%%%%%%%%%%%%
where
%%%%%%%%%%%%%%%%%%%
\begin{equation}
\label{eq66}
\tau_{m}^{-}= 4m+4m^{-}-u +1
+\frac{2i}{\pi}\ln\left\vert
2\cos\left(\frac{1}{2}u \pi\right)\right\vert.
\end{equation}
%%%%%%%%%%%%%%%%%%%

We can numerically check our asymptotic approximations for the zeros by plugging them into the identity (see (\ref{eq29b}))
%%%%%%%%%%%%%%%%%%%
\begin{equation}
\label{eq66a}
1+\frac{e^{\frac13(3u-1)\pi i}\mathrm{Ai}_{1}
\left(\mathpzc{a}_{m}\right)}
{\mathrm{Ai}_{-1}\left(\mathpzc{a}_{m}\right)}=0.
\end{equation}
%%%%%%%%%%%%%%%%%%%

If (for $m$ not sufficiently large) the accuracy is not satisfactory, we can easily improve it by numerically solving for small complex $\delta_{m}$ in the identity
%%%%%%%%%%%%%%%%%%%
\begin{equation}
\label{eq66d}
1+\frac{e^{\frac13(3u-1)\pi i}\mathrm{Ai}_{1}
\left(\mathpzc{a}_{m}^{\mathrm{approx}}+\delta_{m}\right)}
{\mathrm{Ai}_{-1}
\left(\mathpzc{a}_{m}^{\mathrm{approx}}+\delta_{m}\right)}=0,
\end{equation}
%%%%%%%%%%%%%%%%%%%
where $\mathpzc{a}_{m}^{\mathrm{approx}}$ was our original approximation via the asymptotic expansions given above, giving $\mathpzc{a}_{m}^{\mathrm{approx}}+\delta_{m}$ as our improved approximation. We find this to be very effective.

\section*{Acknowledgement}
The authors acknowledge support from {\emph{Ministerio de \allowbreak Econom\'{\i}a y Competitividad}}, project PID2021-127252NB-I00 (MCIN/AEI/10.13039    \allowbreak  
/501100011033/ FEDER, UE.

\bibliographystyle{siamplain}
\bibliography{biblio}

\end{document}